\documentclass[12pt]{article}
\usepackage[small,compact]{titlesec}
\usepackage{amsmath,amssymb,amsfonts,mathabx,setspace}
\usepackage{graphicx,caption,epsfig,subfigure,epsfig,wrapfig}
\usepackage{url,color,verbatim,algorithmicx}
\usepackage[noend]{algpseudocode}
\usepackage[ruled,boxed]{algorithm} 

\usepackage{enumitem} 
\usepackage{booktabs}  
\usepackage[scaled]{helvet}
\usepackage[T1]{fontenc}
\usepackage[bookmarks=true, bookmarksnumbered=true, colorlinks=true,   pdfstartview=FitV,
linkcolor=blue, citecolor=blue, urlcolor=blue]{hyperref}
\usepackage[english]{babel}  
\usepackage{multirow}
\usepackage[numbers,sort]{natbib} 
\usepackage[normalem]{ulem}
\usepackage{mathtools}
\newcommand{\mathsout}[1]
{\bgroup\mathchoice
  {\sbox0{$\displaystyle{#1}$}%
    \usebox0\hspace{-\wd0}%
    \rule[0.5\ht0-0.5\dp0-.5pt]{\wd0}{1pt}}%
  {\sbox0{$\textstyle{#1}$}%
    \usebox0\hspace{-\wd0}%
    \rule[0.5\ht0-0.5\dp0-.5pt]{\wd0}{1pt}}%
  {\sbox0{$\scriptstyle{#1}$}%
    \usebox0\hspace{-\wd0}%
    \rule[0.5\ht0-0.5\dp0-.5pt]{\wd0}{1pt}}%
  {\sbox0{$\scriptscriptstyle{#1}$}%
    \usebox0\hspace{-\wd0}%
    \rule[0.5\ht0-0.5\dp0-.5pt]{\wd0}{1pt}}%
\egroup}
\def\yolambda{\lambda}

\def\by{\mathbf{y}} 
 
 \def\bw{\mathbf{w}} 
\def\bL{\mathbf{L}}
\def\bA{\mathbf{A}}  \def\bC{\mathbf{C}}
\def\bB{\mathbf{B}}\def\bb{\mathbf{b}}
\def\be{\mathbf{e}}

\def\R{\mathbb{R}}

\def\costF{\mathcal{E}}

\newcommand{\Ebracket}[1]{\mathbb{E}\left[{#1}\right]}
\def\E{\mathbb{E}}

\newcommand{\norm}[1]{\left\|#1\right\|}

\newcommand{\innerp}[1]{\langle{#1}\rangle}

\newcommand{\mathspan}[1]{ \mathrm{span}\left\{ {#1} \right\} }

\newcommand{\argmin}[1]{\underset{#1}{\operatorname{arg}\operatorname{min}}\;}

\def\calE{\mathcal{E}}
\def\mathspan{\mathrm{span}}
\def\Gbar{ {\overline{G}} }

\def\LGbar{ {\mathcal{L}_{\overline{G}}}  }

\def\calS{\mathcal{S}}
\def\calT{\mathcal{T}}

\def\spaceY{{L^2_\mu(\mathcal{T})}} 

\newcommand{\beqa}{\begin{eqnarray}}
\newcommand{\eeqa}{\end{eqnarray}}
\newcommand{\beqas}{\begin{eqnarray*}}
\newcommand{\eeqas}{\end{eqnarray*}}
\newcommand{\beq}{\begin{equation}}
\newcommand{\eeq}{\end{equation}}
\newcommand{\beqs}{\begin{equation*}}
\newcommand{\eeqs}{\end{equation*}}
\newcommand{\one}{{\rm{\mathbf{1}}}}
\newtheorem{theorem}{Theorem}

\newtheorem{definition}[theorem]{Definition}

\newtheorem{lemma}[theorem]{Lemma}
\newtheorem{proposition}[theorem]{Proposition}
\newtheorem{remark}[theorem]{Remark}
\newenvironment{proof}[1][Proof]{\noindent\textbf{#1.} }{\ \rule{0.5em}{0.5em}}

\numberwithin{equation}{section}
\numberwithin{theorem}{section}

\usepackage{xcolor}
\definecolor{darkmagenta}{rgb}{0.55, 0.0, 0.55}
\colorlet{colorYO}{darkmagenta}

\newcommand{\bphi}{\boldsymbol{\phi}}
\newcommand{\bpsi}{\boldsymbol{\psi}}
\begin{document}

\begin{center}
{\Large 
An adaptive RKHS regularization for Fredholm integral equations
} \\[0pt] 
\vspace{4mm}
Fei Lu and Miao-Jung Yvonne Ou
\footnote{FL: Department of Mathematics, Johns Hopkins University; feilu@math.jhu.edu  \\
MYO: Department of Mathematical Science, University of Delaware; mou@udel.edu }

\end{center}
 
\begin{abstract} 
 Regularization is a long-standing challenge for ill-posed linear inverse problems, and a prototype is the Fredholm integral equation of the first kind with additive Gaussian measurement noise. 
We introduce a new RKHS regularization adaptive to measurement data and the underlying linear operator. This RKHS arises naturally in a variational approach, and its closure is the function space in which we can identify the true solution. Also, we introduce a small noise analysis to compare regularization norms by sharp convergence rates in the small noise limit. Our analysis shows that the RKHS- and $L^2$-regularizers yield the same convergence rate when their optimal hyper-parameters are selected using the true solution, and the RKHS-regularizer has a smaller multiplicative constant. However, in computational practice, the RKHS regularizer 
significantly outperforms the $L^2$-and $l^2$-regularizers in producing consistently converging estimators when the noise level decays or the observation mesh refines.    
\end{abstract}


\section{Introduction}

We consider the inverse problem of recovering the input function in the Fredholm integral equation of the first kind from \emph{discrete noisy} output. Specifically, let $\calS\subset \R$ and $\calT\subset \R$ be two compact sets. For a given integral kernel $K$, 
we aim to recover the function $\phi:\calS\to \R$ from a discrete noisy dataset $\by = \{y(t_i),  t_0<t_1 <\cdots <  t_m\}$, which is generated from 
\begin{align}\label{eq:FIE}
y(t) = \int_\calS K(t,s)\phi(s) \nu(ds)+ \sigma \dot W(t) =: L_K\phi(t) + \sigma \dot W(t). 
\end{align}
Here $\nu$ is a finite measure on $\calS$, and it is suitable for both continuous and discrete models: it is the Lebesgue measure when $\calS$ is an interval and an atom measure when $\calS$ has finitely many elements. Similarly, we define $\mu$ as an atomic measure with $\mu(t_i) =t_{i+1}-t_i$, which can be viewed as the discrete approximation of the Lebesgue measure. 
The measurement noise  $\sigma \dot W(t)$ is the white noise; that is, the noise at $t_i$ has a Gaussian distribution $\mathcal{N}(0, \sigma^2 (t_{i+1}-t_i))$ for each $i$. Such a noise is integrable when the observation mesh refines, i.e., $\max_i(t_{i+1}-t_i)$ vanishes. 

The integral operator $L_K$ in \eqref{eq:FIE} maps from $L^2_\nu(\calS)$ to $L^2_\mu(\calT)$.  The given integral kernel $K: \calT\times \calS \to \R$ is assumed to be continuous. A typical example is $ K(t,s) = s^{-2}e^{- s t}$
 with $\calS \subset [a,b]$ and $\calT\subset [0,T]$ for some $b>a>0$ and $T>0$, which arises from the magnetic resonance relaxometry (MRR) \cite{Bi2022span-of-regular}; see Section \ref{sec:num} for more details.

Eq.\eqref{eq:FIE} is a prototype of ill-posed inverse problems, dating back to Hadamard \cite{hadamard1923lectures}, and it remains a testbed for regularization methods \cite{wahba1973convergence,engl1996regularization,hansen1998rank,Li2005modified,gazzola2019ir}. It is ill-posed in $L^2_\nu(\calS)$ in the sense that the least squares solution with mini-norm is sensitive to measurement noise in data \cite{nashed1974generalized}. In other words,  in the variational formulation of the inverse problem, the \emph{minimal norm minimizer} of a loss functional
\beqa
&&\mathcal{E}(\phi):=\norm{L_K\phi-\by}^2_{\spaceY} 
\label{least_sq_org}
\eeqa
 is sensitive to the data $\by$. Such an ill-posedness roots in the fact that the operator $L_K: L^2_\nu(\calS)\to L^2_\mu(\calT)$ is compact with eigenvalues converging to zero.

 Regularization, necessary to avoid amplifying the measurement noise, is a long-standing challenge for ill-posed inverse problems, and numerous regularization methods exist. Roughly speaking, in the variational formulation, a regularization method restricts the space of search either by constraints on $\phi$ or the loss functional (e.g., minimizing $\calE(\phi)$ subject to $\|\phi\|_\square\leq \alpha_0$ or minimizing $\|\phi\|_\square$ subject to $\calE(\phi)\leq \alpha_1$); or by adding a penalty to the loss functional $\mathcal{E}$, e.g., the well-known Tikhonov regulation:  
\beqa\label{eq:regu_loss}
\phi_{\yolambda} = \argmin{\phi} \calE_{\yolambda}(\phi ) \mbox{ with } \calE_{\yolambda}(\phi ) : =  \mathcal{E}(\phi)+{\yolambda} \|\phi\|_\square^2. 
\eeqa 
Here $\|\phi\|_\square$ is a certain regularization norm and $(\alpha_0,\alpha_1,{\yolambda})$ are the hyper-parameters. There is tremendous effort in selecting the norm and the hyper-parameters. Once the norm is chosen, the minimization can be solved by either direct or iterative methods. The direct methods solve the linear equations by canonical matrix decompositions, and the selection of ${\yolambda}$ has been thoroughly studied (see e.g.,\cite{wahba1977practical,hansen1998rank,hansen_LcurveIts_a}). In particular, the L-curve method in \cite{hansen_LcurveIts_a} is a widely used robust method. The iterative methods are flexible for high-dimensional problems, and we refer to \cite{gazzola2019ir,chen2022stochastic,zhang2022stochastic} and the reference therein for recent developments. 

In this paper, we focus on the selection of the regularization norm. In particular, we introduce a new data-adaptive Reproducing Kernel Hilbert Space (DA-RKHS) norm and compare it with the widely used $L^2$-norm for regularization. 


The new RKHS-regularization is inspired by the data-adaptive RKHS regularization for learning the kernels in operators \cite{LLA22,LAY22,LangLu21,chada2022data},
 and arises naturally in the variational approach. The closure of this RKHS is referred to as the \emph{function space of identifiability} (FSOI), where a unique solution can be identified. Specifically, we introduce an exploration measure $\rho$ on $\calS$ to quantify how the kernel explores
 the function $\phi(s)$ at different $s$. The variational inverse problem is then formulated in $L^2_\rho(\calS)$. 

The RKHS norm will result in a specific regularization norm $\|\phi\|_\square$ that restricts the minimization search inside the RKHS, preventing errors outside the FSOI from contaminating the regularized solution.

Also, we present \emph{sharp convergence rates} for the RKHS- and $L^2$-regularized estimators in the small noise limit.
These convergence rates are sharp in the sense that they come with explicit multiplicative constants; see Theorem \ref{thm:conv_sigma}. These sharp estimates enable a theoretical comparison between the two regularizers: the RKHS- and $L^2$-regularized estimators both converge linearly, but the RKHS-regularizer outperforms the $L^2$-regularizer with a smaller multiplicative factor. 

Furthermore, we implement the new RKHS regularization in a matrix-decomposition-based algorithm (see Sec.\ref{sec:num}), which is a variant of DARTR in \cite{LLA22}. This algorithm uses the L-curve method in \cite{hansen_LcurveIts_a} to select the optimal hyper-parameter in the Tikhonov regularization.
 Numerical results show that the RKHS-regularizer consistently outperforms the $l^2$- and $L^2$-regularizer when either the noise level decays or the observation mesh $\max_i(t_{i+1}-t_i)$ refines. In particular, the RKHS norm results in better L-shaped curves than the other two norms, indicating that the RKHS norm provides a better metric on the function space for the success of the L-curve method.
 The MATLAB code for our numerical tests is available at \url{https://github.com/LearnDynamics/dartr_Fredholm.git}.

There are various regularization norms, including the Euclidean norms (e.g., \cite{hansen1998rank,tihonov1963solution}), the total variation norm $\|\phi'\|_{L^1}$ in Rudin--Osher--Fatemi method in \cite{rudin1992nonlinear}, the $L^1$ norm $ \|\phi\|_{L^1}$ in LASSO (e.g., \cite{tibshirani1996_RegressionShrinkage}), and the RKHS-norm with a user-specified reproducing kernel  (e.g., \cite{wahba1977practical,bauer2007regularization,CZ07book}) with tuned hyper-parameters for the reproducing kernel. However, these norms are based on pre-assumed properties of the solution and are generic without considering the specific inverse problem. In contrast, our RKHS norm is adaptive to the operator and the observation mesh; our reproducing kernel is intrinsic to the inverse problem and it is implicitly determined without any hyper-parameter to be tuned.  


\subsection{Main contributions}
This study has two main contributions. 
        \begin{itemize}
        \item \emph{Introducing a new adaptive RKHS regularization}. We introduce a new adaptive RKHS for Tikhonov regularization. This RKHS has a reproducing kernel determined by the integral kernel and data, and its closure is the linear subspace in which the inverse problem has a unique solution. Thus, the resulting RKHS regularization is adaptive. We propose an implementation algorithm based on matrix decomposition. Also, we numerically demonstrate that the RKHS regularization outperforms the commonly-used $L^2$ and $l^2$ regularization by yielding more accurate estimators with more consistent convergence when either the noise level decays or the observation mesh refines.              
       \item \emph{Introducing a small noise analysis approach for comparing estimators using different regularization norms}. We introduce a small noise analysis to establish a theoretical comparison between two regularization norms: the RKHS and $L^2$ norms. The analysis shows that when the optimal $\lambda$ is selected by minimizing the expected error, the convergence rates of the two are the same. However, the adaptive RKHS regularizer outperforms the $L^2$ regularizer with a smaller multiplicative constant.
        \end{itemize}
        
         To the best of our knowledge, this study is the first to use the adaptive RKHS for regularizing Fredholm-type inverse problems, and it is the first result comparing regularization norms through sharp convergence rates in the small noise limit.
         The adaptive RKHS (with a reproducing kernel determined by the linear model and data) for regularization applies to general linear machine learning and inverse problems. In particular, it provides an automatic reproducing kernel for the Gaussian process and kernel-based regression. Also, the comparative analysis in this study sets a theoretical framework for comparing regularization norms, which will lead to a better understanding of regularization and the selection of regularization norms.  
         The limitations of this study are discussed in Section \ref{sec:limitations}. 

\bigskip
The rest of the paper is organized as follows. We introduce the adaptive RKHS in Section \ref{sec:FSOI}, with a characterization of the RKHS and the function space of identifiability. Section \ref{sec:conv} proves the convergence of the RKHS-regularized estimator, and Section \ref{sec:num} presents the algorithm and numerically demonstrates the robust convergence of the estimator when either the noise decays or the observation mesh refines. We discuss the limitations of this study in Section \ref{sec:limitations} and conclude with future developments of the adaptive RKHS regularization strategy in Section \ref{sec:conlusion}.   

\section{An adaptive RKHS for the inverse problem} \label{sec:FSOI}

In this section, we introduce the RKHS based on the variational formulation of the inverse problem. 
A core element is an identifiability theory specifying the function space in which the loss functional has a unique minimizer. Importantly, if the true solution is in this space, the minimizer of the loss functional recovers the solution when there is no measurement noise in the data. Thus, we call this \emph{function space of identifiability} (FSOI). For this FSOI, an RKHS norm is introduced and used as the regularization norm $\| \phi\|_\square$ in \eqref{eq:regu_loss}. We will show that this choice of penalty term ensures the search of the minimum of \eqref{eq:regu_loss} takes place in the FSOI.

\subsection{The function space of identifiability} 
We first introduce an ambient function space $L^2_\rho(\calS)$, where the measure $\rho$ is defined as  
\beq
\label{exp_measure}
\frac{d\rho}{d\nu}(s):=\frac{1}{Z} \int_\calT | K(t,s)| \mu(dt), \, \forall s\in \calS, 
\eeq
where $Z= \int_\calS \int_\calT | K(t,s) |\mu(dt)\nu(ds)$ is the normalizing constant. From \eqref{eq:FIE}, we see that $\rho$ quantifies the exploration by the integral kernel $K$ to the unknown input function through the output set $\calT$, and hence is referred to as an \emph{exploration measure}.  

The major advantage of the space $L^2_\rho(\calS)$ over the original space $L^2_\nu(\calS)$ is that it is adaptive to the specific setting of the inverse problem. In particular, this weighted space takes into account the structure of the integral kernel and the data points in $\calT$. It can provide better scaling and reduce the ill-conditioning in computation (see \cite{LangLu21}). Thus, while the following identifiability theory can be carried out for both $L^2_\rho(\calS)$ and $L^2_\nu(\calS)$, we will focus only on $L^2_\rho(\calS)$.

We define the function space of identifiability by the loss functional.
 \begin{definition} The function space of identifiability (FSOI) is the largest linear subspace of $L^2_\rho(\calS)$ where the loss functional $\mathcal{E}$ in \eqref{least_sq_org} has a unique minimizer. 
 \end{definition}
 
The FSOI is the space in which the Fr\'echet derivative of the quadratic loss functional has a unique zero. This motivates us to study the Fr\'{e}chet derivative of the loss functional  $\mathcal{E}$ in $L^2_\rho(\calS)$ and identify the operator of inversion when solving the zero of the derivative of $\mathcal{E}$. We start with writing the quadratic term $\innerp{L_K\phi,L_K\phi}_{\spaceY}$ in the loss functional \eqref{least_sq_org} into a bilinear form
\beqa 
\label{double_innerp}
\innerp{L_K\phi,L_K\psi}_{\spaceY}= \int_\calS \int_\calS \phi(s)\psi(s')G(s,s')\nu(ds) \nu(ds'),
\eeqa
where the integral kernel $G:\calS\times \calS \to \R$ is defined as
 \beq \label{eq:G}
  G(s,s'):= \int_\calT K(t,s)K(t,s') \mu(dt).  \eeq
Note that 
\begin{equation}\label{eq:bd_G}
G(s,s')\leq \|K\|_\infty Z \min\{\frac{d\rho}{d\nu}(s),\frac{d\rho}{d\nu}(s')\},\, \forall(s,s')\in \calS\times \calS.
\end{equation}
 Then, we can weigh it by the exploration measure and define 
\beq
\Gbar(s,s'):=
 \frac{G(s,s')}{\frac{d\rho}{d\nu}(s)\frac{d\rho}{d\nu}(s')}. 
\label{Gbar_def}
\eeq

The continuity of the kernel $K\in C(\calT\times \calS)$ implies the following lemma.
\begin{lemma}
\label{mercer_kernel_prop} 
Assume that $K\in C(\calT\times \calS)$ and recall $\rho$ in \eqref{exp_measure}. The following statements hold.
\begin{itemize}
\item[{\rm(a)}] The function $\Gbar$ in  \eqref{Gbar_def} is a Mercer kernel. The operator $\LGbar: L^2_\rho(\calS) \to L^2_\rho(\calS)$ defined by 
\beq
\label{LGbar_def}
\LGbar\phi(r):=\int_a^b \phi(s)\Gbar(r,s)\rho(ds)
\eeq
is compact, self-adjoint and positive. Moreover, for $\phi,\psi \in L^2_\rho(\calS)$, we have 
\beq
\innerp{L_K\phi,L_K\psi}_{\spaceY} = \innerp{\LGbar \phi,\psi}_{L^2_\rho(\calS)}=\innerp{ \phi,\LGbar \psi}_{L^2_\rho(\calS)}.
\label{quadratic_2}
\eeq
\item[{\rm(b)}] Let $\{\lambda_i\}_{i\geq 1}$ be the positive eigenvalues of $\LGbar$ arranged in a descending order. Denote by $\{\psi_i,\,\psi_j^0\}_{i,j}$ the orthonormal eigenfunctions of $\LGbar$ with $\psi_i$ and $\psi_j^0$ corresponding to eigenvalues $\lambda_i$ and zero (if any), respectively. Then, $\{\lambda_i\}_{i\geq 1}$ is either finite or $\lambda_i\to 0$ as $i\to \infty$, and these eigenfunctions form a complete basis of $L^2_\rho(\calS)$. 
\item[{\rm(c)}] The operator $\LGbar$ is a trace-class operator,  
\begin{equation}\label{eq:trace}
\sum_i \lambda_i = \int_\calS \Gbar(s,s)\rho(ds) \leq  \|K\|_\infty Z \nu(\calS)<\infty, 
\end{equation}
where $\nu(\calS)$ is the measure of $\calS$ and $Z$ is the normalizing constant in \eqref{exp_measure}. 
\end{itemize}
\end{lemma}

 \begin{proof} 
Since the kernel $K$ is continuous, so is $\rho$ and $G$ in \eqref{eq:G}, thus $\Gbar$ is also continuous.  Note that $\Gbar$ is symmetric, i.e., $\Gbar(s,s')=\Gbar(s', s)$ for any $s,s'\in \calS$. Also, it is positive semi-definite, i.e., 
 \beq \sum_{i=1}^n\sum_{j=1}^n c_i c_j \Gbar(s_i,s_j) =\int_\calT \left(\sum_{i=1}^n  c_i \one_{\calS}(s_i)  K(t,s_i)/\frac{d\rho}{d\nu}(s_i) \right)^2 \, \mu(dt) \ge 0
 \eeq
 for any $\{c_j\}_{j=1}^n\subset \R$, $\{s_j\}_{j=1}^k\subset \calS$ and $n\in \mathbb{N}$. 
 Moreover, by \eqref{eq:bd_G}, we have $\Gbar\in L^2(\rho\times\rho)$ because  
  $$\int_\calS \int_\calS \Gbar(s,s') ^2\rho(ds) \rho(ds') \leq \|K\|_\infty^2 Z^2 <\infty.
  $$ Thus, $\LGbar$ is compact and positive self-adjoint (e.g., \cite[Proposition 4.6]{CZ07book}).

Part (b) follows directly from that $\LGbar$ is compact and positive self-adjoint (\cite[Theorem 4.7]{CZ07book}).

To prove Part (c), recall that from Mercer's theorem, we have 
\[
\Gbar(s,s')=\sum_i \lambda_i \psi_i(s)\psi_i(s'),
\]
where the convergence is uniform on $\calS\times \calS$. Then, 
 \begin{align*}
 \sum_i\lambda_i = \int_\calS \Gbar(s,s)\rho(ds) = \int_\calS \frac{G(s,s)}{\frac{d\rho}{d\nu}(s)}ds \leq \|K\|_\infty Z \nu(\calS), 
 \end{align*}
 where the inequality follows from \eqref{eq:bd_G}. 
 \end{proof} 

The next theorem characterizes the FSOI through the Fr\'echet derivative of the loss functional, highlighting the role of $\LGbar$ as the operator of inversion in the variational approach. 
\begin{theorem}[Function space of identifiability]\label{thm:FSOI}
Assume that $K\in C(\calT\times \calS)$ and let $\by$ be the data generated with a function $\phi_{true}\in L^2_\rho(\calS)$ by \eqref{eq:FIE}. Then, we have the following characterization of the inverse problem. 
\begin{itemize}
\item[{\rm(a)}]  There exists a unique $\phi^\by \in L_\rho^2(\calS)$ such that
\beq
\innerp{\phi^\by,\psi}_{L^2_\rho(\calS)}=\innerp{L_K\psi,\by}_{\spaceY},\quad \forall \psi \in L^2_\rho(\calS), \label{linear}
\eeq
and it has a decomposition 
\begin{equation}\label{eq:noise_dec}
\phi^\by = \LGbar \phi_{true}+ \phi^\sigma,
\end{equation} 
where $\phi^\sigma$ has a distribution $\mathcal{N}(0, \sigma^2\LGbar)$, i.e., $\phi^\sigma= \sum_i \sigma \xi_i \lambda_i^{1/2}\psi_i $ with $\{\xi_i\}$ being independent identically distributed standard Gaussian random variables and $\Ebracket{\|\phi^\sigma\|_{L^2_\rho(\calS)}^2}= \sigma^2\sum_i\lambda_i$. 
\item[{\rm(b)}] The Fr\'echet derivative of the loss functional in $L^2_\rho$ is $\nabla \calE(\phi) = 2(\LGbar \phi- \phi^\by)$, where $\phi^\by$ is defined in \eqref{linear}. 
\item[{\rm(c)}] The FSOI of $\calE$ is $H:=\overline{ \mathspan\{\psi_i\}_{i:\lambda_i>0} }$ with the closure taken in $L^2_\rho(\calS)$. 
\item[{\rm(d)}] Assume that the data is noisy, i.e. $\sigma> 0$. If $\sum_{i:\lambda_i>0} \lambda_i^{-1}<\infty$, the unique minimizer of $\calE$ in $H$ is $\widehat{\phi} = \LGbar^{-1}\phi^\by$; but if $\sum_{i:\lambda_i>0} \lambda_i^{-1}=\infty$, the solution $\LGbar^{-1}\phi^\by$ is ill-defined in $L^2_\rho(\calS)$  in the sense that $\Ebracket{\|\LGbar^{-1}\phi^\by\|_{L^2_\rho(\calS)}^2} = \infty$. 
\item[{\rm(e)}]  When the observation is noiseless, we have $\widehat{\phi} = \LGbar^{-1}\phi^\by= P_H\phi_{true}$, where $P_H$ is the projection operator of $H$. 
\end{itemize}
\end{theorem}

\begin{proof}
The existence and uniqueness of $\phi^\by$ in Part $(a)$ follow from the Riesz representation theorem. Denote the observation by $\by= L_K\phi_{true}+\sigma \Delta \mathbf{W}$ with $\Delta \mathbf{W} = \left(W(t_{i})- W(t_{i-1}): 1\leq i\leq m \right)$. Then, the decomposition of $\phi^\by$ follows from \eqref{linear}: 
\begin{align*}
\innerp{\phi^\by,\psi}_{L^2_\rho(\calS)} = \innerp{L_K\psi,L_K\phi_{true}}_\spaceY +   \innerp{L_K\psi,\sigma \Delta \mathbf{W}}_\spaceY 
& = \innerp{\LGbar\phi_{true},\psi}_{L^2_\rho(\calS)} + \innerp{\psi,\phi^\sigma}_{L^2_\rho(\calS)} ,
\end{align*}
where the first term in the last equation comes from the definitions of the operator $\LGbar$ in \eqref{LGbar_def}, and the second term comes from the Riesz representation. The distribution of $\phi^\sigma$ is $\mathcal{N}(0,\sigma^2\LGbar)$ in the sense that for each $\psi \in L^2_\rho(\calS)$, the random variable $ \innerp{\psi,\phi^\sigma}_{L^2_\rho(\calS)} =   \innerp{L_K\psi,\sigma \Delta \mathbf{W}}_\spaceY  $ is Gaussian with mean zero and variance 
$$
{\E[ \innerp{\psi,\phi^\sigma}_{L^2_\rho(\calS)}^2] = \E[\innerp{L_K\psi,\sigma \Delta \mathbf{W}}_\spaceY^2  ] = \sigma^2\E[\innerp{L_K\psi,L_K\psi}_\spaceY  ] = \sigma^2\innerp{\psi,\LGbar\psi}_{L^2_\rho(\calS)}.
}$$
 Therefore, we can write $\phi^\sigma= \sum_i \sigma \xi_i \lambda_i^{1/2}\psi_i $ with $\{\xi_i\}$ being i.i.d.~standard Gaussian, and $\Ebracket{\|\phi^\sigma\|_{L^2_\rho(\calS)}^2}= \sigma^2\sum_i\lambda_i$, where the sum is finite by \eqref{eq:trace}.

Part $(b)$ follows directly from the definition of the Fr\'echet derivative. In fact, by \eqref{quadratic_2} and \eqref{linear}, we can write the loss functional as 
\begin{equation}\label{eq:loss_L2}
\calE(\phi) =  \innerp{\LGbar \phi,\phi}_{L^2_\rho(\calS)} +  \innerp{\phi^\by,\phi}_{L^2_\rho(\calS)} + \|\by\|_\spaceY^2. 
\end{equation}
Then, the Fr\'echet derivative $\nabla \calE(\phi)$ in $L^2_\rho(\calS)$ is 
\[
\innerp{\nabla \calE(\phi), \psi}_{L^2_\rho(\calS)} = \lim_{h\to 0} \frac{\calE(\phi+h\psi) - \calE(\phi)}{h} =  \innerp{2 ( \LGbar \phi -\phi^\by), \psi}_{L^2_\rho(\calS)}, \, \forall \psi \in L^2_\rho(\calS). 
\]

For Part $(c)$, first, note that the quadratic loss functional has a unique minimizer in $H$ because its derivative has a unique zero in $H$. Meanwhile, note that $H$ is the orthogonal complement of the null space of $\LGbar$, and $\calE(\phi_{true}+\phi^0) = \calE(\phi_{true})$ for any $\phi^0$ such that $\LGbar \phi^0=0$.  Thus, $H$ is the largest such linear subspace of $L^2_\rho(\calS)$, and we conclude that  $H$ is the FSOI.
 
For Part (d), note that 
\[  \LGbar^{-1}\phi^\sigma = \sum_{i:\lambda_i>0} \lambda_i^{-1/2} \sigma \xi_i \psi_i .
\]
Then,  if $\sum_{i:\lambda_i>0} \lambda_i^{-1}<\infty$, which happens only when there are finitely many non-zero eigenvalues, we have $ \LGbar^{-1}\phi^\sigma \in L^2_\rho(\calS)$ because $\Ebracket{ \| \LGbar^{-1}\phi^\sigma\|_{L^2_\rho(\calS)}^2} =\sigma^2  \sum_{i:\lambda_i>0} \lambda_i^{-1}<\infty$ {by using the fact that $\{\xi_i\}$ are i.i.d.~ standard Gaussian.} Thus, $\phi^\sigma$ in $\LGbar(L^2_\rho)$, so is $\phi^\by$, and the estimator $\widehat{\phi} = \LGbar^{-1}\phi^\by$ is well-defined in $L^2_\rho$. By Part $(b)$, this estimator is the unique zero of the loss functional's Fr\'echet derivative in $H$. Hence it is the unique minimizer of $\calE(\phi)$ in $H$. 

In contrary, if $\sum_{i} \lambda_i^{-1}=\infty$, we have $\Ebracket{\| \LGbar^{-1}\phi^\sigma \|_{L^2_\rho(\calS)}^2 }
 = \sum_{i} \lambda_i^{-1} \sigma^2 
=\infty$, and hence $\LGbar^{-1}\phi^\by$ is ill-defined.

For Part (e), when the data is noiseless, i.e., $\by=L_K\phi_{true}$, we have $\phi^\by= \LGbar\phi_{true}$ from Part $(a)$. Hence $\widehat{\phi} = \LGbar^{-1}\phi^\by = P_H \phi_{true}$. 
\end{proof}

Theorem \ref{thm:FSOI} reveals the nature of the ill-posedness of this inverse problem, and provides insights on regularization.  The variational inverse problem is ill-defined outside the FSOI $H$. Its ill-posedness in $H$ depends on the smallest eigenvalue of the operator $\LGbar$.  
\begin{itemize}
\item When the data is noiseless, the minimizer of the loss functional is the $H$-projection of the true input function. In other words, the inverse problem can only recover the $H$-projection of the true input function. 
\item When the data is noisy, its minimizer $\LGbar^{-1}\phi^\by$ can be ill-defined in $L^2_\rho$ if $\sum_{i:\lambda_i>0} \lambda_i^{-1}=\infty$, but it is well-defined when $\sum_{i:\lambda_i>0} \lambda_i^{-1}<\infty$.
\end{itemize}
Based on these observations, the regularization scheme proposed in Section \ref{sec:RKHS} is designed to ensure the solution lies in the FSOI and either remove the non-integrable components related to small eigenvalues or reduce the bias caused by the noise.   

\subsection{An adaptive RKHS regularization}\label{sec:RKHS}

We introduce an adaptive RKHS Tikhonov regularization. This RKHS is defined based on the operator of inversion $\LGbar$; thus, it is adaptive to the original integral operator $L_K$ and the data settings.  The next lemma characterizes this RKHS.  

Recall that for a given Mercer Kernel $G:X\times X \rightarrow \mathbb{R}$ on a metric space $X$, $G_x:=G(x,\cdot)$, a Hilbert space $(H_G,(\cdot,\cdot)_{H_G})$ is an RKHS of functions on $X$ with reproducing kernel $G$ if it has all the following three properties: (i) for all $x\in X$, $G_x\in H_G$, (2) the span of $\{G_x, \, x\in X\}$ is dense in $H_G$, (iii) for all $f\in H_G$ and $x\in X$, $f(x)=(G_x,f)_{H_G}$, cf. \cite{cucker2007learning-theory}.

\begin{lemma}[Characterization of the adaptive RKHS]\label{lemma:rkhs} 
Assume $K\in C(\calT\times\calS)$. The RKHS $H_G$  with $\Gbar$ as its reproducing kernel satisfies the following properties. 
\begin{enumerate}
\item[{\rm(a)}]  $H_G:=\LGbar^{\frac{1}{2}}(L^2_\rho(\calS))$ and its inner product is 
\begin{equation*} 
\innerp{\phi,\psi}_{H_G}:= \innerp{\LGbar^{-\frac{1}{2}}\phi,\LGbar^{-\frac{1}{2}}\psi}_{L^2_\rho(\calS)}. 
\end{equation*}
The operator $\LGbar$ is self-adjoint in $H_G$. Moreover, we have
$\innerp{\phi,\psi}_{L^2_\rho(\calS)}=\innerp{\LGbar \phi,\psi}_{H_G}$  for any $\phi\in L^2_\rho(\calS)$ and $\psi\in H_G$. 
\item [{\rm(b)}] 
  $\{\sqrt{\lambda_i} \psi_i\}_{i=1}^\infty$ is an orthonormal basis of $H_G$, where $\{\psi_i\}_{i}$ are the orthonormal eigenfunctions of $\LGbar$ in $L^2_\rho(\calS)$ corresponding to positive eigenvalues $\{\lambda_i\}_{i}$. 
\item [{\rm(c)}] 
For any $\phi=\sum_{i=1}^\infty c_i \psi_i$, with $c_i\in \R$, we have 
\beq
\innerp{L_K\phi,L_K\phi}_\spaceY=\sum_{i=1}^\infty \lambda_i c_i^2,\quad \| \phi\|^2_{L^2_\rho}=\sum_{i=1}^\infty c_i^2, \quad \| \phi\|^2_{H_G}=\sum_{i=1}^\infty \lambda_i^{-1} c_i^2 \text{ if } \phi\in H_G. 
\label{norms}
\eeq  
Moreover, the $H_G$ norm is stronger than the $L^2_\rho$-norm: $\norm{\phi}^2_{H_G}\ge \lambda_1^{-1} \norm{\phi}^2_{L^2_\rho}$.
\item [{\rm(d)}]  $H=\overline{ H_G}$ with closure in $L^2_\rho(\calS)$, where $H=\overline{ \mathspan\{\psi_i\}_{i:\lambda_i>0} }$ is the FSOI. 
\end{enumerate}
\end{lemma}
 \begin{proof}
 The first part is from the standard characterization theorem of RKHS, e.g., \cite[Section 4.4]{CZ07book} or \cite{LLA22}. When $\phi\in L^2_\rho(\calS)$, we have $\LGbar\phi \in H_G$. Then, by the definition of the inner product and the symmetry of $\LGbar^{-1/2}$, we have $\innerp{\LGbar\phi,\psi}_{H_G} = \innerp{\LGbar^{1/2}\phi,\LGbar^{-1/2}\psi}_{L^2_\rho(\calS)} = \innerp{\phi,\psi}_{L^2_\rho(\calS)}$. 

Part \rm{(b)} follows directly from the characterization of the inner product in Part (a). Part (c) follows directly from \eqref{quadratic_2}, the orthonormality of the eigenfunctions and the characterization of the inner product. Part \rm{(d)} {follows from that $\overline{H_G}=\overline{ \mathspan\{\sqrt{\lambda_i}\psi_i\}_{i:\lambda_i>0} }= \overline{ \mathspan\{\psi_i\}_{i:\lambda_i>0} }= H$. }
 \end{proof}

Next, we propose to use the RKHS norm for regularization, that is, use $\| \phi\|_{H_G}$ in the place of $\| \phi\|_{\square}$ in \eqref{eq:regu_loss}. Note that by using the RKHS norm, the loss functional is minimized over $H_G$, not $H\subset L^2_\rho(\calS)$, because elements in $H$ may have unbounded $H_G$ norm.  

\begin{proposition}[Regularized estimators] \label{prop_error} Consider the $L^2_\sigma$-regularized estimator $\widehat{\phi}_{{\alpha}}^{L^2_\rho}$, and the $H_G$-regularized estimator $\widehat{\phi}_{{\alpha}}^{H_G}$, which can be expressed in terms of $\LGbar$ as follows: 
\begin{align}
\widehat{\phi}_{\yolambda}^{L^2_\rho}&:= \argmin{\phi\in L^2_\rho} \mathcal{E}(\phi)+{\yolambda} \norm{\phi}_{L^2_\rho}^2 = (\LGbar+{\yolambda} I)^{-1}\phi^\by,\label{regs_L2}\\ 
\widehat{\phi}_{\yolambda}^{H_G}&:= \argmin{\phi\in H_G} \mathcal{E}(\phi)+{\yolambda} \norm{\phi}_{H_G}^2 = (\LGbar^2+{\yolambda} I)^{-1}\LGbar \phi^\by. 
 \label{regs_HG}
\end{align}
Express the true function in terms of the eigen-basis functions of $\LGbar$ as $\phi_{true}= \sum_i c_i\psi_i + \sum_j d_j\psi_j^0$.   Then, for any ${\yolambda}>0$, the biases of these two regularized estimators satisfy the following estimates:  
 \beqa
  && \norm{\widehat{\phi}_{\yolambda}^{L^2_\rho}-\phi_{true}}^2_{L^2_\rho} =   \sum_i (\lambda_i +{\yolambda})^{-2}(\sigma \lambda_i^{1/2}\xi_i - {\yolambda} c_i)^2+ \sum_jd_j^2,\label{errL2_eig}\\ 
 && \norm{\widehat{\phi}_{\yolambda}^{H_G}-\phi_{true}}^2_{L^2_\rho}  =   \sum_i (\lambda_i^2 +{\yolambda})^{-2}(\sigma \lambda_i^{3/2}\xi_i - {\yolambda} c_i)^2+ \sum_jd_j^2, \label{errHG_eig}
  \eeqa
  where $\xi_i$ are i.i.d.~standard Gaussian random variables in Theorem {\rm \ref{thm:FSOI}}.  
\end{proposition}
\begin{proof}
The minimizers' uniqueness and explicit form follow from the Fr\'echet derivatives of the regularized loss functionals. In fact, by Theorem \ref{thm:FSOI}(b),  the Fr\'echet derivative of $\calE_{\yolambda}(\phi)  = \mathcal{E}(\phi)+{\yolambda} \norm{\phi}_{L^2_\rho(\calS)}^2$ in $L^2_\rho(\calS)$ is  
\[
\nabla \calE_{\yolambda}(\phi) = 2[ (\LGbar + {\yolambda} I) \phi -\phi^\by ]. 
\]
It has a unique zero,  $\widehat{\phi}_{\yolambda}^{L^2_\rho} = (\LGbar + {\yolambda} I)^{-1}\phi^\by$, which is the unique minimizer of $\calE_{\yolambda}$ in $L^2_\rho$, and thus we have \eqref{regs_L2}. Similarly, note that by Lemma \ref{lemma:rkhs} {\rm(a)} and \eqref{eq:loss_L2},  for any $\phi\in H_G$, 
\[
\calE(\phi)  
=  \innerp{\LGbar^2\phi,\phi}_{H_G} +  \innerp{\LGbar\phi^\by,\phi}_{H_G} +  \|\by\|_\spaceY^2. 
\]
Then, the Fr\'echet derivative of $\widetilde \calE_{\yolambda}  = \mathcal{E}(\phi)+{\yolambda} \norm{\phi}_{H_G}^2$ in $H_G$ is  
\[
\nabla \widetilde \calE_{\yolambda}(\phi) = 2[ (\LGbar^2 + {\yolambda} I) \phi -\LGbar\phi^\by ]. 
\]
Its unique zero,  $\widehat{\phi}_{\yolambda}^{H_G} = (\LGbar^2 + {\yolambda} I)^{-1}\LGbar\phi^\by$, is the minimizer of $\widetilde \calE_{\yolambda}$ in $H_G$ as in \eqref{regs_HG}. 

The eigenvalue characterizations of the biases in \eqref{errL2_eig} and \eqref{errHG_eig} follow from the decomposition of $\phi^\by = \LGbar \phi_{true}+\phi^\sigma$ with $\phi^\sigma= \sum_i \sigma \lambda_i^{1/2}\xi_i \psi_i$ in \eqref{eq:noise_dec} and the facts listed below: 
\begin{align*}
& \widehat{\phi}_{\yolambda}^{L^2_\rho}=(\LGbar+{\yolambda} I)^{-1}(\LGbar{\phi_{true}}+{\yolambda} \phi_{true}-{\yolambda} \phi_{true}+\phi^\sigma)=\phi_{true}+(\LGbar+{\yolambda} I)^{-1}(-{\yolambda} \phi_{true}+\phi^\sigma),\\
& \widehat{\phi}_{\yolambda}^{H_G}=(\LGbar^2+{\yolambda} I)^{-1}(\LGbar^2{\phi_{true}}+\LGbar\phi_1^\delta)=\phi_{true}+(\LGbar^2+{\yolambda} I)^{-1}(-{\yolambda} \phi_{true}+\LGbar\phi^\sigma).
\end{align*}
Then, we obtain \eqref{errL2_eig} and \eqref{errHG_eig}.  
\end{proof}

Proposition \ref{prop_error} demonstrates the complexity of choosing an optimal hyper-parameter ${\yolambda}$. An optimal ${\yolambda}$ aims to balance the error caused by the noise and the shift from the true solution caused by the regularization, and it depends on the spectrum of the operator, each realization of the noise, and the true solution. Thus, it is important to select an optimal ${\yolambda}$ adaptive to these factors, and the L-curve method does so.  

Also, Proposition \ref{prop_error} shows that the true solution's components outside the FSOI remain a bias for both estimators, because there is no information about these components in the data. Thus, in the error analysis presented in Section \ref{sec:conv}, we focus only on the error inside the FSOI. 

\section{Convergence of the regularized estimators as the noise decays}
\label{sec:conv}

A key question in regularization is the choice of the regularization norm $\| \phi \|_\square$ in \eqref{eq:regu_loss} when there is little prior knowledge. Unfortunately, empirical comparisons do not lead to a conclusive answer due to multiple mixed factors, including the random noise, the spectrum of the inversion operator, the smoothness of the solution, the model/discretization error, the selection of the hyper-parameter, and the method of optimization (matrix decomposition for least squares or iterative methods for the minimization). Therefore, one must extract the essentials from these multiple factors to make a meaningful comparison of different regularized estimators.

         We introduce next a comparative analysis framework to compare regularization norms in the small noise limit. The small noise limit approach helps us separate the aforementioned mixed factors and analyze their effects in regularization. Hence, this comparative analysis framework enables us to compare the adaptive RKHS norm and the $L^2_\rho$ norm, making a small step toward guiding the regularization norm selection in practice.

We show that both the RKHS- and $L^2_\rho$-regularized estimators converge in mean-squares at order  $1$ as $\sigma\downarrow 0$ under the assumption that the spectrum of $\LGbar$ converges exponentially or polynomially. The RKHS-regularizer slightly outperforms the $L^2_\rho$-regularizer by a smaller multiplicative factor. The spectral assumption and the convergence are observed in the numerical tests in Section \ref{sec:num}. 

In the following comparative analysis, we assume $\phi_{true}= \sum_i c_i\psi_i\in H$. 
We estimate the mean-square error (MSE) of the $L^2_\rho$-regularized and the $H_G$-regularized estimators:
\[
e^{L^2_\rho}({\yolambda}) := \E \norm{\widehat{\phi}_{\yolambda}^{L^2_\rho}-\phi_{true}}^2_{L^2_\rho},\quad e^{H_G}({\yolambda}) :=  \E \norm{\widehat{\phi}_{\yolambda}^{H_G}-\phi_{true}}^2_{L^2_\rho}. 
\]
Note that by Proposition \ref{prop_error},  these MSEs are 
 \begin{equation}\label{eq:mse} 
 \begin{aligned}
 e^{L^2_\rho}({\yolambda})  &=   \sum_i (\lambda_i +{\yolambda})^{-2}(\sigma^2 \lambda_i + {\yolambda}^2 c_i^2),\\ 
 e^{H_G}({\yolambda})&=   \sum_i (\lambda_i^2 +{\yolambda})^{-2}(\sigma^2 \lambda_i^{3}+ {\yolambda}^2 c_i^2). 
 \end{aligned}
   \end{equation}
The purpose of the following theorem is to give a quantitative comparison of these two regularized estimators.   

\begin{theorem}[Small noise limits of regularized estimators]
\label{thm:conv_sigma}
 {Consider the inverse problem \eqref{eq:FIE} with a true solution }
 $\phi_{true}= \sum_i c_i\psi_i \in H= \overline{\mathrm{span}\{\psi_i\}}$, where $\{\psi_i\}$ are the eigen-functions of positive eigenvalues of $\LGbar$ in \eqref{LGbar_def}. 
If the decay rate of the spectrum of $\LGbar$, which is determined by the integral kernel $K$, is either exponential or polynomial, 
 the following upper bounds and sharp rates hold for the MSEs in \eqref{eq:mse} as the standard deviation of the Gaussian noise $\sigma\to 0$. 
\begin{itemize}
\item \textbf{Exponential spectrum decay.}  Assume that the spectrum of the operator $\LGbar$ decays to 0 exponentially, i.e., there exists $\theta>0$ such that $\lambda_i = e^{-\theta i}$ for all $i\geq 1$. 
\begin{enumerate}
	\item [(a)]
	If $\phi_{true}$ satisfies the condition $M_0:=\sup_{i\geq 1}\lambda_i^{-1}c_i^2 <\infty$, which is valid for every function in $H_G$ but can also hold for functions not belonging to $H_G$, then we have the following estimates: 
	\begin{equation*}
	\begin{aligned}
	\min_{{\yolambda}>0}e^{H_G}({\yolambda}) &\leq e^{H_G}(\sigma^2) \leq (1+ M_0) \frac{\pi}{4\theta} \sigma+ O(\sigma^2);\\ 
	\min_{{\yolambda}>0}e^{L^2_\rho}({\yolambda}) &\leq e^{L^2_\rho}(\sigma) \leq (1+ M_0) \frac{2}{\theta} \sigma+ O(\sigma^2), 
	\end{aligned}
	\end{equation*} 
	where $O(\sigma^2)$ is the big-O notation. 
	\item[(b)] 
	
	 In particular, for $\phi_{true}$ satisfying $c_i^2 =\lambda_i$, we have \textbf{sharp rates} of convergence: 
	\begin{equation*}
	\begin{aligned}
	 e^{H_G}({\yolambda}_{opt}) 
	&=& \frac{\pi}{4\theta} \sigma+ O(\sigma^2) \quad \text{ with }{\yolambda}_{opt}&:= \argmin{{\yolambda}>0} e^{H_G}({\yolambda}) = \sigma^2, \\ 
	 e^{L^2_\rho}(\widetilde {\yolambda}_{opt}) 
	&=& \frac{2}{\theta}\sigma  + O(\sigma^2) \quad \text{ with } \widetilde {\yolambda}_{opt}&:= \argmin{{\yolambda}>0} e^{L^2_\rho}({\yolambda}) = \sigma + O(\sigma^2). 
	\end{aligned}
	\end{equation*} 
\end{enumerate}

\item \textbf{Power spectrum decay.} Assume that the spectrum of the operator $\LGbar$ decays polynomially, i.e., there exists $\theta>0$ such that $\lambda_i = i^{-\theta}$ for all $i\geq 1$.
	\begin{enumerate}
	\item[(c)] 
	 If $M_0:=\sup_{i\geq 1}\lambda_i^{-1}c_i^2 <\infty$, then we have the following estimates:
	 \begin{equation*}
		\begin{aligned}
		\min_{{\yolambda}>0}e^{H_G}({\yolambda}) &\leq e^{H_G}(\sigma^2) \leq (1+ M_0)\sigma \left( C_{H_G}+ o(1) \right);\\ 
		\min_{{\yolambda}>0}e^{L^2_\rho}({\yolambda}) &\leq e^{L^2_\rho}(\sigma) \leq (1+ M_0)\sigma (C_{L^2_\rho} + o(1))), 
		\end{aligned}
		\end{equation*}
		where $C_{H_G}$ and $C_{L^2_\rho}$ are two positive constants. 
	\item[(d)] In particular, for $\phi_{true}$ satisfying $c_i^2 =\lambda_i$, 
 we have \textbf{sharp rates} of convergence:
		\begin{equation*}
		\begin{aligned}
		e^{H_G}({\yolambda}_{opt}) & =& \sigma \left( C_{H_G}+ o(1) \right) \quad \text{ with } 	{\yolambda}_{opt}&= \argmin{{\yolambda}>0} e^{H_G}({\yolambda}) = \sigma^2,  \\ 
		e^{L^2_\rho}(\widetilde {\yolambda}_{opt})  & =& \sigma (C_{L^2_\rho} + o(1)) \quad \text{ with } 	\widetilde {\yolambda}_{opt}& = \argmin{{\yolambda}>0} e^{L^2_\rho}({\yolambda}) =\sigma (  C_{\widetilde {\yolambda}_{opt}}+ o(1)),
		\end{aligned}
		\end{equation*}
		where the constants are {	$C_{H_G} =\frac{1}{2}C_{L^2_\rho} =\frac{1}{2\theta}  \Gamma(\frac{1}{2}-\frac{1}{2\theta}) \Gamma(\frac{1}{2}+\frac{1}{2\theta})$ and $C_{\widetilde {\yolambda}_{opt}}=\frac{\sqrt{\theta+1}}{\sqrt{\theta-1}}$.}  
	\end{enumerate}
\end{itemize}

\end{theorem}

\begin{remark}[Sharp rates and restrictive assumptions] The sharp rates in Parts {(b)} and {(d)} rely on the restrictive assumptions on the coefficients of $\phi_{true}$. These sharp rates show that the RKHS regularizer outperforms the $L^2_\rho$ regularizer in the sense that it has a smaller multiplicative constant. On the other hand, these restrictive assumptions can be significantly relaxed when a sharp rate is not needed. For example, one only needs a uniform bound for the Picard ratio in Parts {\rm (a)} and {\rm (c)}  when one only aims for an upper bound; one can also allow perturbations in the form of  $\lambda_i= p_i i^{-\theta}$ or $\lambda_i= p_i e^{-\theta i}$ with {the perturbations $0<a \leq p_i\leq b<\infty$ for all $i$ for some $a,b>0$} when one is only concerned with the rate but not the multiplicative constant  {\rm(see \cite{LangLu23sna} for more details)}. 
\end{remark}

To prove the theorem, we first introduce some notations and two technical lemmas. We rewrite the MSEs in \eqref{eq:mse} as 
\begin{equation}\label{eq:error_lambda}
e^{H_G}({\yolambda}) = \sigma^2 A({\yolambda}) + {\yolambda}^2 B({\yolambda}); \quad e^{L^2_\rho}({\yolambda}) = \sigma^2 \widetilde A({\yolambda}) + {\yolambda}^2 \widetilde B({\yolambda}). 
\end{equation}
where $A$ and $B$ are defined by
\begin{equation}\label{eq:AB}
\begin{aligned}
A({\yolambda}) &= \sum_i (\lambda_i^2+{\yolambda})^{-2} \lambda_i^{3}, \quad & B({\yolambda}) &= \sum_i (\lambda_i^2+{\yolambda})^{-2}  c_i^2, \\
\widetilde A({\yolambda}) &= \sum_i (\lambda_i+{\yolambda})^{-2} \lambda_i, \quad &\widetilde B({\yolambda}) & = \sum_i (\lambda_i+{\yolambda})^{-2} c_i^2. 
\end{aligned}
\end{equation}
Then, since a minimizer must be a critical point, we have that ${\yolambda}_{opt}$ must satisfy
\[
0=\frac{d}{d{\yolambda}} e^{H_G}({\yolambda})  =\sigma^2 A'({\yolambda}) + 2{\yolambda} [ B({\yolambda}) +\frac{{\yolambda}}{2}B'({\yolambda})]
\]
and similarly for $\widetilde {\yolambda}_{opt}$. With $B_1({\yolambda}) : =B({\yolambda}) +\frac{{\yolambda}}{2} B'({\yolambda}) $ and $\widetilde B_1({\yolambda}) : = \widetilde B({\yolambda}) +\frac{{\yolambda}}{2} \widetilde B'({\yolambda}) $, we obtain the following lemma providing equations for the minimizers of the MSEs. 
\begin{lemma}\label{lemma:AB-lambda}
The minimizers of $e^{H_G}$ and $e^{L^2_\rho}$ in \eqref{eq:mse}, denoted by ${\yolambda}_{opt}$ and $\widetilde {\yolambda}_{opt}$, respectively, satisfy
\begin{equation}\label{eq:lambda_opt}
{\yolambda}_{opt} =- \sigma^2\frac{A'({\yolambda}_{opt})}{2B_1({\yolambda}_{opt})}, \quad \widetilde {\yolambda}_{opt} =- \sigma^2\frac{\widetilde A'(\widetilde {\yolambda}_{opt})}{2\widetilde B_1(\widetilde {\yolambda}_{opt})}, 
\end{equation}
where the functions $A'({\yolambda})$ and $A'({\yolambda})$ $B_1({\yolambda})$, $\widetilde A({\yolambda})$ and $\widetilde B_1({\yolambda})$ are given by 
\begin{equation}\label{eq:A1B1}
\begin{aligned}
A'({\yolambda}) &= -2\sum_i (\lambda_i^2+{\yolambda})^{-3} \lambda_i^{3}, \quad &B_1({\yolambda}) =  \sum_i (\lambda_i^2+{\yolambda})^{-3} \lambda_i^2 c_i^2, \\
\widetilde A'({\yolambda}) &= -2\sum_i (\lambda_i+{\yolambda})^{-3} \lambda_i, \quad &\widetilde B_1({\yolambda})=  \sum_i (\lambda_i+{\yolambda})^{-3}\lambda_i c_i^2. \\
\end{aligned}
\end{equation}
\end{lemma}

The next lemma estimates these series by Riemann sum when ${\yolambda}$ is small when the spectrum decays exponentially. 
\begin{lemma}\label{lemma:seiresAB}
Assume that $\lambda_i = e^{-\theta i}$ for all $i\geq 1$ with $\theta>0$. Then, for small ${\yolambda}>0$, we have 
\begin{equation}\label{eq:seriesAB}
\begin{aligned}
A({\yolambda})=  \sum_i (\lambda_i^2+{\yolambda})^{-2} \lambda_i^{3} & = \frac{1}{2\theta \sqrt{{\yolambda}}} [\arctan\frac{1}{\sqrt{{\yolambda}}} - \frac{\sqrt{{\yolambda}}}{1+{\yolambda}}] + O(1), \\ 
 B_c({\yolambda}): = \sum_i (\lambda_i^2+{\yolambda})^{-2} \lambda_i & = \frac{1}{2\theta }{\yolambda}^{-3/2} [\arctan\frac{1}{\sqrt{{\yolambda}}} + \frac{\sqrt{{\yolambda}}}{1+{\yolambda}}] + O(1), \\ 
\widetilde A({\yolambda})=   \sum_i (\lambda_i+{\yolambda})^{-2} \lambda_i & = \frac{1}{\theta {\yolambda}(1+{\yolambda})}+ O(1), \\ 
\widetilde A'({\yolambda})=  -2    \sum_i (\lambda_i+{\yolambda})^{-3} \lambda_i & = \frac{-(1+2{\yolambda})}{\theta {\yolambda}^2(1+{\yolambda})^2}+ O(1), \\ 
\widetilde   B_c({\yolambda}):     =   \sum_i (\lambda_i+{\yolambda})^{-2} \lambda_i^2 & = \frac{1}{2\theta {\yolambda} (1+{\yolambda})^2}+ O(1). \\ 
\end{aligned}
\end{equation}
\end{lemma}
\begin{proof}
The proof is based on the Riemann sum approximation of integrals. Note that for $k\in \{1,3\}$, the function $f(x) = ( e^{-2\theta x} + {\yolambda})^{-2} e^{-\theta xk } $ satisfies $\int_0^\infty f(x) dx= \sum_{i=1} f(i)  + O(1)$. Thus, 
\begin{align*}
 & \sum_i (\lambda_i^2+{\yolambda})^{-2} \lambda_i^k  = \int_0^\infty  ( e^{-2\theta x} + {\yolambda})^{-2} e^{-\theta xk }  dx + O(1) =  \frac{1}{\theta} \int_0^1 \frac{t^{k-1} dt }{(t^2+{\yolambda})^2}  + O(1), 
\end{align*}
where we applied a change of variables $t= e^{-\theta x}$ to obtain the second equality. 
Then, the first two equations in \eqref{eq:seriesAB} follow directly from the facts that $\int_0^1 \frac{t^2 dt }{(t^2+{\yolambda})^2} = \frac{1}{2 \sqrt{{\yolambda}}} [\arctan\frac{1}{\sqrt{{\yolambda}}} - \frac{\sqrt{{\yolambda}}}{1+{\yolambda}}] $ and $\int_0^1 \frac{dt }{(t^2+{\yolambda})^2} = \frac{1}{2}{\yolambda}^{-3/2} [\arctan\frac{1}{\sqrt{{\yolambda}}} + \frac{\sqrt{{\yolambda}}}{1+{\yolambda}}] $. 

Similarly, we obtain the last three equations in \eqref{eq:seriesAB} by using the integrals $\int_0^1 \frac{1}{(t+{\yolambda})^2}dt $,  $\int_0^1 \frac{1}{(t+{\yolambda})^3}dt $ and $\int_0^1 \frac{t}{(t+{\yolambda})^3}dt $. 
\end{proof}

\vspace{2mm}

The next lemma estimates these series when the spectrum decays polynomially. 
\begin{lemma}\label{lemma:seiresAB_power}
Assume that $\lambda_i = i^{-\theta}$ for all $i\geq 1$ with $\theta>0$. Denote $C_\theta(s,k,\alpha) = \Gamma(\gamma)\Gamma(k - \gamma) $ with $\gamma= \frac{1}{1+s}(\alpha-\frac{1}{\theta})$, where $s\in \{0,1\}$, $\alpha\in \{1,2,3\}$, and $k\in \{2,3 \}$. Then, for small ${\yolambda}>0$, we have 
\begin{equation}\label{eq:seriesAB_power}
\begin{aligned}
A({\yolambda})=  \sum_i (\lambda_i^2+{\yolambda})^{-2} \lambda_i^{3} & =  \frac 1 {2\theta} {\yolambda}^{\gamma-2
} [C_\theta(1,2,3)+ o(1)] \, \text{with } \gamma  = \frac{1}{2}(3-\frac{1}{\theta}); 
\\ 
 B_c({\yolambda}): = \sum_i (\lambda_i^2+{\yolambda})^{-2} \lambda_i & = \frac 1 {2\theta} {\yolambda}^{\gamma-2
}  [C_\theta(1,2,1)+ o(1)] \, \text{with } \gamma  = \frac{1}{2}(1-\frac{1}{\theta}); \\ 
\widetilde A({\yolambda})=   \sum_i (\lambda_i+{\yolambda})^{-2} \lambda_i & =\frac 1 {\theta} {\yolambda}^{\gamma-2
}  [C_\theta(0,2,1)+ o(1)]\, \text{with } \gamma  = (1-\frac{1}{\theta});  \\ 
\widetilde A'({\yolambda})=  -2    \sum_i (\lambda_i+{\yolambda})^{-3} \lambda_i & = -\frac 1 {\theta} {\yolambda}^{\gamma-3
}  [C_\theta(0,3,1)+ o(1)]\, \text{with } \gamma  = (1-\frac{1}{\theta});  \\ 
\widetilde   B_c({\yolambda}):     =   \sum_i (\lambda_i+{\yolambda})^{-3} \lambda_i^2 & = \frac 1 {2\theta} {\yolambda}^{\gamma-3
}  [C_\theta(0,3,2)+ o(1)] \, \text{with } \gamma  = (2-\frac{1}{\theta}). \\ 
\end{aligned}
\end{equation}
\end{lemma}
\begin{proof} Again, the proof is based on the Riemann sum approximation of integrals. 

To estimate the five terms in unified notations, let 
\[
F_{s}({\yolambda},\alpha,k) :=  \sum_i (\lambda_i^{1+s}+{\yolambda})^{-k} \lambda_i^{\alpha}, \quad \gamma= \frac{\alpha-\frac{1}{\theta}}{1+s}.
\]
Then, with $s=1$, we have $A({\yolambda}) = F_{1}({\yolambda},3,2)$ and $B_c({\yolambda}) =  F_{1}({\yolambda},1,2)$, which use $k=2$ and $\alpha \in\{1,3\}$. 
Also, with $s=0$, we have $\widetilde A({\yolambda}) = F_0({\yolambda},1,2)$, $\widetilde A'({\yolambda}) =  F_0({\yolambda},1,3)$, and  $\widetilde B_c({\yolambda}) =  F_0({\yolambda},2,3)$ using $\alpha\in\{1,2\}$ and $k\in \{2,3\}$. 

Note that the function $f(x) = ( x^{-(1+s)\theta} + {\yolambda})^{-k} x^{-\theta \alpha } $ satisfies $\int_1^\infty f(x) dx= \sum_{i=1} f(i)  + O(1)$. Thus, under the assumption $\lambda_i=i^{-\theta}$, we have
\begin{align*}
 F_{s}({\yolambda},\alpha,k) &=  \sum_i (\lambda_i^{1+s}+{\yolambda})^{-k} \lambda_i^\alpha  = \int_1^\infty  ( x^{-\theta(1+s)} + {\yolambda})^{-k} x^{-\theta \alpha }  dx + O(1) \\
 &=  \frac{1}{\theta} \int_0^1 \frac{t^{\alpha-1 -1/\theta} dt }{(t^{1+s}+{\yolambda})^k}  + O(1), 
\end{align*}
where we applied a change of variables $t= x^{-\theta}$ (equivalently, $x=t^{-1/\theta}$ and $dx=-\frac{1}{\theta} t^{-1-1/\theta}dt$) to obtain the second equality. We will estimate the last integral with $s=1$ and $s=0$ separately.  

Consider first the cases with $s=1$, i.e., the series $A({\yolambda}) = F_{1}({\yolambda},3,2)$ and $B_c({\yolambda}) =  F_{1}({\yolambda},1,2)$, which use $k=2$ and $\alpha \in\{1,3\}$. 
By a change of variables $z= t^2/{\yolambda}$ (and hence $t= {\yolambda}^{1/2}  z^{1/2}$ and $dt =\frac 1 2 {\yolambda}^{1/2} z^{-1/2} dz$), and note that $\gamma=\frac{1}{2}(\alpha-\frac{1}{\theta})<2 $, we have 
\begin{align*}
 F_{1}({\yolambda},\alpha,2)& =\frac{1}{\theta} \int_0^1 \frac{t^{\alpha-1 -1/\theta} dt }{(t^2+{\yolambda})^2} = \frac 1 {2\theta} {\yolambda}^{\gamma-2} \int_0^{1/{\yolambda}} (z+1)^{-2}z^{\gamma-1}dz \\
 &= \frac 1 {2\theta} {\yolambda}^{\gamma-2} \left(\frac{\Gamma(\gamma)\Gamma(2 - \gamma)}{\Gamma(2)} + o(1)\right),
\end{align*}
where the last equality follows from that fact that $ \int_0^{\infty} (z+1)^{-k}z^{\gamma-1}dz= \frac{\Gamma(\gamma)\Gamma(k - \gamma)}{\Gamma(k)}$ and 
$\int_0^{1/{\yolambda}} (z+1)^{-2}z^{\gamma-1}dz = \int_0^{\infty} (z+1)^{-2}z^{\gamma-1}dz +o(1)$.  
Recall that $\Gamma(2)= 1$. Then, the first two equations in \eqref{eq:seriesAB_power} follow.  

Next, consider the cases with $s=0$, i.e., the series $\widetilde A({\yolambda}) = F_0({\yolambda},1,2)$, $\widetilde A'({\yolambda}) =  F_0({\yolambda},1,3)$, and  $\widetilde B_c({\yolambda}) =  F_0({\yolambda},2,3)$ using $\alpha\in\{1,2\}$ and $k\in \{2,3\}$. By a change of variables $z=t/{\yolambda}$, and note that $\gamma=\alpha-\frac{1}{\theta}<k$,  we have
\begin{align*}
 F_{0}({\yolambda},\alpha,k) &=\frac{1}{\theta} \int_0^1 \frac{t^{\alpha-1 -1/\theta} dt }{(t+{\yolambda})^k} =  \frac{1}{\theta}{\yolambda}^{\gamma-k} \int_0^{1/{\yolambda}} (z+1)^{-k}z^{\gamma-1}dz \\
&= \frac 1 \theta {\yolambda}^{\gamma-k} \left(\frac{\Gamma(\gamma)\Gamma(k - \gamma)}{\Gamma(k)} + o(1)\right). 
\end{align*}
Then,  we obtain the last three equations in \eqref{eq:seriesAB_power} using the facts that $\Gamma(2)=1$ and $\Gamma(3)=2$. 
\end{proof}

\vspace{3mm}
\begin{proof}[Proof of Theorem \ref{thm:conv_sigma}]
We prove Parts (a) and (c) by solving algebraic equations to obtain the optimal hyper-parameters using Lemmas \ref{lemma:seiresAB}-\ref{lemma:seiresAB_power}, and we prove Parts (b) and (d) by direct evaluations.

\textbf{Part (a).} The upper bounds follow from direct evaluations of $ e^{H_G}(\sigma^2)$ and $ e^{L^2_\rho}(\sigma)$. In fact, recall that $e^{H_G}(\sigma^2)= \left[ \sigma^2A({\yolambda}) + {\yolambda}^2B({\yolambda}) \right] \mid_{{\yolambda} = \sigma^2}$, and the bound for  $\sigma^2A(\sigma^2) $ comes from  \eqref{eq:seriesAB}:
\[
\sigma^2 A(\sigma^2) = \sigma^2 \left( \frac{1}{2\theta \sqrt{{\yolambda}}} [\arctan\frac{1}{\sqrt{{\yolambda}}} - \frac{\sqrt{{\yolambda}}}{1+{\yolambda}}] + O(1) \right)  \mid_{{\yolambda} = \sigma^2} \leq \frac{\pi}{4\theta} \sigma +o(\sigma).  
\]
Also, recalling that  $M_0= \sup_{i} \lambda_i^{-1}c_i^2 <\infty$, wa have, 
\begin{align*}
{\yolambda}^2 B({\yolambda})  \mid_{{\yolambda} = \sigma^2} & = {\yolambda}^2 \sum_i (\lambda_i^2+{\yolambda})^{-2}  \lambda_i  \lambda_i^{-1}c_i^2 \mid_{{\yolambda} = \sigma^2} \leq M_0  {\yolambda}^2 \sum_i (\lambda_i^2+{\yolambda})^{-2}  \lambda_i  \mid_{{\yolambda} = \sigma^2}\\
  & \leq  \frac{M_0}{2\theta } [ {\yolambda}^{1/2} [\arctan\frac{1}{\sqrt{{\yolambda}}} + \frac{\sqrt{{\yolambda}}}{1+{\yolambda}}] + O(1) ] \mid_{{\yolambda} = \sigma^2} \leq \frac{M_0\pi}{4\theta} \sigma +o(\sigma). 
 \end{align*}
  where the last inequality follows from the second equation in \eqref{eq:seriesAB}.
 Thus, 
  \begin{align*}
\min_{{\yolambda}>0} e^{H_G}({\yolambda}) \leq  e^{H_G}(\sigma^2) = \left[ \sigma^2A({\yolambda}) + {\yolambda}^2B({\yolambda}) \right] \mid_{{\yolambda} = \sigma^2}  \leq (1+ M_0)  \frac{\pi}{4\theta}\sigma+ O(\sigma^2).
 \end{align*}
 
 Similarly, recall that $e^{L^2_\rho}(\sigma)=  \sigma^2 \widetilde A(\sigma) + \sigma^2\widetilde B(\sigma) $ with $\widetilde A$ and $\widetilde B$ defined in \eqref{eq:AB}. Note that 
 $$\widetilde B({\yolambda}) =\sum_i (\lambda_i+{\yolambda})^{-2} c_i^2 = \sum_i (\lambda_i+{\yolambda})^{-2} \lambda_i \lambda_i^{-1}  c_i^2  \leq \sum_i (\lambda_i+{\yolambda})^{-2} \lambda_i M_0 = M_0 \widetilde A({\yolambda}). 
 $$
Meanwhile, the upper bound for $\widetilde A$ in \eqref{eq:seriesAB} implies that  
\[ \sigma^2 \widetilde A(\sigma)\leq \sigma^2\left[ \frac{1}{\theta\sigma (1+\sigma)} + O(1) \right] \leq \frac{1}{\theta}\sigma + O(\sigma^2).
\]
Hence, 
$e^{L^2_\rho}(\sigma) \leq (1+M_0) \sigma^2 \widetilde A(\sigma) \leq  (1+ M_0) \frac{2}{\theta} \sigma+ O(\sigma^2)$. 

\textbf{Part (b).} Substituting $c_i^2= \lambda_i$ into the equation for $B_1({\yolambda})$ in \eqref{eq:A1B1}, we have $B_1({\yolambda}) =  \sum_i (\lambda_i^2+{\yolambda})^{-3} \lambda_i^3 = \frac{-1}{2}A'({\yolambda})$. Hence, Eq. \eqref{eq:lambda_opt} implies that ${\yolambda}_{opt}= \sigma^2$. Also, substituting  $c_i^2= \lambda_i$ into the equation for $B({\yolambda})$ in \eqref{eq:AB}, we obtain
 \[
 B({\yolambda}) = \sum_i (\lambda_i^2+{\yolambda})^{-2}  \lambda_i =  \frac{1}{2\theta }{\yolambda}^{-3/2} [\arctan\frac{1}{\sqrt{{\yolambda}}} + \frac{\sqrt{{\yolambda}}}{1+{\yolambda}}] + O(1). 
 \]
 where the second equality follows from the second equation in \eqref{eq:seriesAB} when ${\yolambda}$ is small (since we will set it to be $\sigma^2$). This equation, together with the first equation in  \eqref{eq:seriesAB} and \eqref{eq:error_lambda}, implies 
 \begin{align}
& \min_{{\yolambda}>0} e^{H_G}({\yolambda}) = e^{H_G}(\sigma^2) = \left[ \sigma^2A({\yolambda}) + {\yolambda}^2B({\yolambda}) \right] \mid_{{\yolambda} = \sigma^2}  \notag \\
 = &  \left[ \frac{ \sigma^2}{2\theta \sqrt{{\yolambda}}} (\arctan\frac{1}{\sqrt{{\yolambda}}} - \frac{\sqrt{{\yolambda}}}{1+{\yolambda}}) + \frac{ 1}{2\theta }\sqrt{{\yolambda}} (\arctan\frac{1}{\sqrt{{\yolambda}}} - \frac{\sqrt{{\yolambda}}}{1+{\yolambda}}) \right] \big|_{{\yolambda}= \sigma^2}+ O(\sigma^2) \notag \\
 = & \sigma^2 \frac{ 1}{2\theta\sqrt{{\yolambda}}  }\arctan\frac{1}{\sqrt{{\yolambda}}} \mid_{{\yolambda}= \sigma^2}+ O(\sigma^2) = \frac{\pi}{4\theta} \sigma+ O(\sigma^2). \notag
 \end{align}
 
 Similarly, to find the minimizer of $e^{L^2_\rho}({\yolambda}) $, we substitute $c_i^2= \lambda_i$ into the equation for $\widetilde B_1({\yolambda})$ in \eqref{eq:A1B1}, we have  $\widetilde B_1({\yolambda}) =  \sum_i (\lambda_i+{\yolambda})^{-3} \lambda_i^2 = \frac{1}{2\theta {\yolambda} (1+{\yolambda})^2}+ O(1)$. Consequently, by the equation for $\widetilde A'({\yolambda})$ in \eqref{eq:seriesAB}, the minimizer $\widetilde {\yolambda}_{opt}$ satisfies  $ {\yolambda}  = -\sigma^2 \frac{\widetilde A'({\yolambda})}{2\widetilde B_1({\yolambda})}$, and it can be solved from 
 \begin{align*}
&&  2 {\yolambda} \widetilde B_1({\yolambda}) & = -\sigma^2 \widetilde A'({\yolambda})   \\
& \Leftrightarrow & \quad   \frac{1}{\theta (1+{\yolambda})^2}+ 2{\yolambda} O(1)&  =  \sigma^2 \frac{1+2{\yolambda}}{\theta {\yolambda}^2 (1+{\yolambda})^2}+  \sigma^2 O(1)   \\ 
& \Leftrightarrow &  {\yolambda}^2 & = \sigma^2+2{\yolambda} \sigma^2+O({\yolambda}^2 \sigma^2 + {\yolambda}^3).  
 \end{align*}
The positive solution to  ${\yolambda}^2= \sigma^2+2{\yolambda} \sigma^2$ is $\sigma^2 + \sigma \sqrt{1+ \sigma^2}$. Thus, the minimizer is $\widetilde {\yolambda}_{opt} =  \sigma + O(\sigma^2)$. Together with the fact that $\widetilde B({\yolambda}) = \sum_i (\lambda_i+{\yolambda})^{-2}\lambda_i = \widetilde A({\yolambda})$, we obtain  
 \begin{align*}
 \min_{{\yolambda}>0} e^{L^2_\rho}({\yolambda}) = &e^{L^2_\rho}({\yolambda}_{opt}) =  \widetilde A({\yolambda}) ( \sigma^2 + {\yolambda}^2 ) \big|_{{\yolambda} = {\yolambda}_{opt}} \\
 = & \frac{1}{\theta {\yolambda} (1+{\yolambda})} ( \sigma^2 + {\yolambda}^2 ) + O(\sigma^2) \big|_{{\yolambda} = {\yolambda}_{opt}} \\
 = & \frac{1}{\theta} \frac{2\sigma^2 + O(\sigma^3)) }{(\sigma +O(\sigma^2)) (1+ \sigma+ O(\sigma^2))} + O(\sigma^2) = \frac{2}{\theta}\sigma + O(\sigma^2), 
 \end{align*}
 where we used Eq.\eqref{eq:seriesAB} to obtain the third equality.

\textbf{Part (c).} The upper bounds in Part (d) follow directly from evaluations of the series using the uniform bound. 

\textbf{Part (d).} The proof is similar to Part (a) with slightly more computations. First, the same argument in Part (a) leads to ${\yolambda}_{opt}= \sigma^2$. To evaluate $e^{H_G}(\sigma^2) = \left[ \sigma^2A({\yolambda}) + {\yolambda}^2B({\yolambda}) \right] \mid_{{\yolambda} = \sigma^2}$, note first that substituting  $c_i^2= \lambda_i$ into the equation for $B({\yolambda})$ in \eqref{eq:AB} and using \eqref{eq:seriesAB_power}, we obtain
 \begin{align*}
 B({\yolambda}) &= \sum_i (\lambda_i^2+{\yolambda})^{-2}  \lambda_i =  B_c({\yolambda}) = \frac{1}{2\theta} {\yolambda}^{-\frac{1}{2}(3+\frac{1}{\theta}) } [C_\theta(1,2,1)+o(1)], \\
 A({\yolambda}) & =  \frac{1}{2\theta} {\yolambda}^{-\frac{1}{2}(1+\frac{1}{\theta}) } [C_\theta(1,2,3)+o(1)]. 
 \end{align*}
 Consequently, we obtain
 \begin{align} \label{eq:errHG_min}
& \min_{{\yolambda}>0} e^{H_G}({\yolambda}) = e^{H_G}(\sigma^2) 
 = \sigma \left( \frac{1}{2\theta}[C_\theta(1,2,1)+ C_\theta(1,2,3)]+ o(1) \right): = \sigma (C_{H_G} + o(1)). 
 \end{align}
Thus, we have $C_{H_G} =\frac{1}{2\theta}[C_\theta(1,2,1)+ C_\theta(1,2,3)] $, which we will compute its value after obtaining a similar expression for $C_{L^2}$.

To find the minimizer of $e^{L^2_\rho}({\yolambda})$, we solve an algebraic equation as in Part (a). 
Substituting $c_i^2= \lambda_i$ into the equation for $\widetilde B_1({\yolambda})$ in \eqref{eq:A1B1}, we have 
 $$\widetilde B_1({\yolambda}) =  \sum_i (\lambda_i+{\yolambda})^{-3} \lambda_i^2 = \widetilde B_c({\yolambda})   = \frac 1 {2\theta} {\yolambda}^{-\frac{1}{\theta}-1}  [C_\theta(0,3,2)+ o(1)]
$$ 
using \eqref{eq:seriesAB_power}. Also, we have $\widetilde A'({\yolambda}) =-\frac 1 {\theta} {\yolambda}^{-2-\frac{1}{\theta}
}  [C_\theta(0,3,1)+ o(1)] $. Hence, since the minimizer $\widetilde {\yolambda}_{opt}$ satisfies  $ {\yolambda}  = -\sigma^2 \frac{\widetilde A'({\yolambda})}{2\widetilde B_1({\yolambda})}$, it can be solved from 
 \begin{align*}
&&  2 {\yolambda} \widetilde B_1({\yolambda}) & = -\sigma^2 \widetilde A'({\yolambda})   \\
& \Leftrightarrow & \quad  {\yolambda}^{-\frac{1}{\theta}}  [C_\theta(0,3,2)+ o(1)]&  =  \sigma^2 {\yolambda}^{-2-\frac{1}{\theta}
}  [C_\theta(0,3,1)+ o(1)] .  
 \end{align*}
Thus, the minimizer is 
\begin{equation}\label{eq:C-lambda_opt}
	\widetilde {\yolambda}_{opt} =  \sigma [ C_\theta(0,3,2)^{-1/2} C_\theta(0,3,1)^{1/2}+ o(1) ]=: \sigma (  C_{\widetilde {\yolambda}_{opt}}+ o(1)). 
\end{equation}
 Together with the fact that $\widetilde B({\yolambda}) = \sum_i (\lambda_i+{\yolambda})^{-2}\lambda_i = \widetilde A({\yolambda}) =\frac 1 {\theta} {\yolambda}^{-1-\frac{1}{\theta}} [C_\theta(1,2,3)+ o(1)] $ as in Eq.\eqref{eq:seriesAB_power}, we obtain  
 \begin{align*}
& \min_{{\yolambda}>0} e^{L^2_\rho}({\yolambda}) = e^{L^2_\rho}({\yolambda}_{opt}) =  \widetilde A({\yolambda}) ( \sigma^2 + {\yolambda}^2 ) \big|_{{\yolambda} = {\yolambda}_{opt}} \\
 = &  \sigma \frac 1 {\theta} C_\theta(1,2,3) [(1+C_\theta(0,3,2)^{-1} C_\theta(0,3,1)) + o(1)] =: \sigma (C_{L^2_\rho} + o(1)). 
 \end{align*} 
 
 At last, we evaluate the constants $C_{H_G}$ in \eqref{eq:errHG_min}, $ C_{\widetilde {\yolambda}_{opt}}$ in \eqref{eq:C-lambda_opt}, and $C_{L^2_\rho}$ above. We have 
 		 \begin{equation*}
		 \begin{aligned}
		 C_{H_G} & = \frac{1}{2\theta}[C_\theta(1,2,1)+ C_\theta(1,2,3)], \\		 C_{\widetilde {\yolambda}_{opt}} & =C_\theta(0,3,2)^{-1/2} C_\theta(0,3,1)^{1/2}, \\
		C_{L^2_\rho}  &= \frac 1 {\theta} C_\theta(1,2,3)[ 1+C_\theta(0,3,2)^{-1} C_\theta(0,3,1)].
		 \end{aligned}
		 \end{equation*}
		Here we denote $C_\theta(s,k,\alpha) = \Gamma(\gamma)\Gamma(k - \gamma) $ with $\gamma= \frac{1}{1+s}(\alpha-\frac{1}{\theta})$, where $s\in \{0,1\}$, $\alpha\in \{ 1+2s, 1+s\}_{s\in \{0,1\}} = \{1,2,3\}$, and $k\in \{2,3 \}$.
				
		Noting that 
		$\gamma=\frac{1}{1+s}(\alpha-\frac 1 \theta) =\begin{cases}
				\frac{1}{2}(3-\frac 1 \theta), & \text{ if } (s,k,\alpha)=(1,2,3); \\
				\frac{1}{2}(1-\frac 1 \theta), & \text{ if } (s,k,\alpha)=(1,2,1),
				\end{cases}
		$ and using the formula 	$\Gamma(z+1) = z\Gamma(z)$, we obtain 
		\begin{align*}
		      C_\theta(1,2,3)&= \Gamma(\frac{3}{2}-\frac{1}{2\theta})\Gamma(\frac{1}{2} + \frac{1}{2\theta}) = \frac{1}{2}(1-\frac{1}{\theta})  \Gamma(\frac{1}{2}-\frac{1}{2\theta}) \Gamma(\frac{1}{2}+\frac{1}{2\theta}), \quad \\
		      C_\theta(1,2,1) &= \Gamma(\frac{1}{2}-\frac{1}{2\theta})\Gamma(\frac{3}{2} + \frac{1}{2\theta}) =  \frac{1}{2}(1+\frac{1}{\theta})  \Gamma(\frac{1}{2}-\frac{1}{2\theta}) \Gamma(\frac{1}{2}+\frac{1}{2\theta}). 	
		\end{align*} 
		Hence, 
		$C_{H_G}  
		 =  \frac{1}{2\theta} C_\theta(1,2,3) [1+C_\theta(1,2,3)^{-1}C_\theta(1,2,1)]
		 = \frac{1}{2\theta}  \Gamma(\frac{1}{2}-\frac{1}{2\theta}) \Gamma(\frac{1}{2}+\frac{1}{2\theta}). 
		$
		
		Similarly, noting that 
		$\gamma=\frac{1}{1+s}(\alpha-\frac 1 \theta) =\begin{cases}
				2-\frac 1 \theta, & \text{ if } (s,k,\alpha)=(0,3,2); \\
				1-\frac 1 \theta, & \text{ if } (s,k,\alpha)=(0,3,1),
				\end{cases},
				$
we have  
				\begin{align*}
		      C_\theta(0,3,2)&= \Gamma(2-\frac{1}{\theta})\Gamma(1+ \frac{1}{\theta}) = (1-\frac{1}{\theta})  \Gamma(1-\frac{1}{\theta}) \Gamma(1+\frac{1}{\theta}), \quad \\
		      C_\theta(0,3,1) &= \Gamma(1-\frac{1}{\theta})\Gamma(2 + \frac{1}{\theta}) = (1+\frac{1}{\theta})  \Gamma(1-\frac{1}{\theta}) \Gamma(1+\frac{1}{\theta}). 
		\end{align*} 
			Consequently, $C_{\widetilde {\yolambda}_{opt}} =C_\theta(0,3,2)^{-1/2} C_\theta(0,3,1)^{1/2}= \frac{\sqrt{\theta + 1}}{\sqrt{\theta-1}}$ and
	\begin{align*}
		 C_{L^2_\rho}  &= \frac 1 {\theta} C_\theta(1,2,3)[ 1+C_\theta(0,3,2)^{-1} C_\theta(0,3,1)]  = \frac{1}{\theta}  \Gamma(\frac{1}{2}-\frac{1}{2\theta}) \Gamma(\frac{1}{2}+\frac{1}{2\theta}). 
	\end{align*}
	Clearly, $C_{H_G} =\frac{1}{2}C_{L^2_\rho}$ and we complete the proof. 
\end{proof}

\section{Algorithm and numerical examples}\label{sec:num}
We present a direct matrix-decomposition-based algorithm to implement the data-adaptive RKHS regularization. We compare the RKHS regularizer with two commonly used regularizers: the $L^2$-regularizer that sets $\|\phi\|_{\square} = \|\phi\|_{L^2_\rho}$ in \eqref{eq:regu_loss} and the $l^2$-regularizer that sets $\|\phi\|_\square = (\sum_{i=1}^n\phi(s_i)^2)^{1/2}$, which is the $l^2$ norm of $\bphi= (\phi(s_1), \ldots,\phi(s_n))$,  the discrete representation of $\phi$ on $\calS$.  

\paragraph{Numerical settings.} 
{
 We consider Eq.\eqref{eq:FIE} with $\calS:= \{s_k\}_{k=1}^n\subset [a,b]$, where $s_k= a+ k\delta $ with $\delta:= (b-a)/n$ for $k\geq 1$ being a uniform mesh. We set $\nu(\{s_k\}) = \delta$ so that \eqref{eq:FIE} can be viewed as a discretization of the following Fredholm integral equation   }  
\begin{align*}
y(t) = \int_{a}^b K(t,s)\phi(s) ds + \sigma \dot W(t)  = :L_K\phi(t) + \sigma \dot W(t), t\in [c,d]. 
\end{align*}
We consider two specific kernels:
\begin{equation}\label{eq:kernels}
\begin{aligned}
K_{exp}(t,s) & = s^{-2}e^{- s t}; \\
K_{poly}(t,s)& = s^{-1}| \sin(s t +1)|,  
\end{aligned}
\end{equation}
where the kernels $K_{exp}$ and $K_{poly}$ lead to operators with exponentially and polynomially decaying spectra, respectively (see Figure \ref{fig:eigVals}). Here the kernel $K_{exp}$ arises from the $T2$ problem in magnetic resonance relaxometry \cite{Bi2022span-of-regular}, where $\phi$ is the distribution of transverse nuclear relaxation times, and it characterizes the material composition of a sample. 

{
We aim to recover $\bphi=  (\phi(s_1), \ldots,\phi(s_n)) \in \R^n$. This discrete vector provides a piecewise-constant approximation to the function $\phi$: $\phi(s)=\sum_{k=1}^n \phi(s_k) \mathbf{1}_{[s_{k-1},s_k]}(s)$ with the convention $s_0=a$, where $\mathbf{1}_A(s)$ is the indicator function of set $A$.  
Then, Eq.\eqref{eq:FIE} can be viewed as a Riemann sum approximation of the above integral: }
\begin{equation}\label{eq:discrete}
\by = \bL \bphi + \bw, 
\end{equation}
where $\bw\in \R^m \sim \mathcal{N}(0,\sigma^2\Delta t I_m)$ and 
\begin{equation}
\bL = \begin{bmatrix}
   K(t_1,s_1)\delta_1 &K(t_1,s_2)\delta_2 & K(t_1,s_3)\delta_3 & \cdots & K(t_1,s_n)\delta_n \\
K(t_2,s_1)\delta_1 &K(t_2,s_2)\delta_2 & K(t_2,s_3)\delta_3 & \cdots & K(t_2,s_n)\delta_n \\
    \vdots & \vdots &  \vdots &  \ddots & \vdots \\
K(t_m,s_1)\delta_1 &K(t_m,s_2)\delta_2 & K(t_m,s_3)\delta_3 & \cdots & K(t_m,s_n)\delta_n \\
    \end{bmatrix}
\end{equation}
with $K(t,s)$ being either $K_{exp}(t,s)$ or $K_{poly}(t,s)$. 

We set $(a,b,c,d)=(1,5,0,5)$ and $n=100$. We set the standard deviation of the noise to be $\sigma= \|\bL \bphi  \|  \times nsr$, where the noise-to-signal ratio $nsr$ is set to be $nsr= 1$ unless otherwise specified. For the data $\by=(y(t_1),\ldots, y(t_m))\in \R^m$, we consider equally spaced grid $ t_i =c+ i\Delta t$, $0\leq i\leq m$ with $\Delta t=(d-c)/m = 0.01$, i.e., with $m=500$, unless otherwise specified.

\paragraph{The space $L^2_\rho(\calS)$.}
The empirical exploration measure of \eqref{eq:discrete}  is 
\begin{equation}\label{eq:rho_example}
\rho(s_k)=\frac{1}{z} \sum_{i=1}^m s_k^{-2}e^{-s_k t_i}\delta, \quad k=1,\ldots,n,
\end{equation}
where $z$ is the normalizing constant. In general, when the true kernel $K$ is available, one may compute the exploration measure analytically and then evaluate $\rho(s_k)$ for each $s_k$ in $\calS$. For example, the analytical exploration measure for $K_{exp}$ is $\rho(s)\, \propto \int_c^d s^{-2}e^{-st}dt = s^{-3}(e^{-cs}- e^{-ds}) $. 

Our goal is to estimate $\bphi\in L^2_\rho(\calS)\subset  \R^n$, which is equipped with the norm $\|f\|_{L^2_\rho(\calS)}^2= \sum_{k=1}^n f_k^2\rho(s_k)$ (and we consider only the values on the support of $\rho$). The support of $\rho$ is where the operator explores the vector $\bphi$, for which we have information for an estimation. Thus, $L^2_\rho(\calS)$ is a proper space with a metric for estimating $\bphi\in \R^n$.

The discrete $\bphi$ uses the Cartesian basis $\{\be_i\}_{i=1}^n$ of  $\R^n$ such that $(\be_1,\be_2,\ldots, \be_n)= I_n$. In $L^2_\rho(\calS)$, these basis vectors have a basis matrix
 \begin{equation}\label{eq:basisMat}
 \bB = ( \innerp{\be_i,\be_j}_{L^2_\rho(\calS)} )_{1\leq i,j\leq n}. 
 \end{equation}
{Note that $\bB=  \mathrm{diag}(\rho(s_k))$ is diagonal for this set of basis vectors, and the inner product of $L^2_\rho(\calS)$ is
$\innerp{\bphi,\boldsymbol{\psi}}_{L^2_\rho}:=\bphi^T \bB \boldsymbol{\psi},\,\text{ for  all } \bphi,\boldsymbol{\psi} \in \R^n$. }

\subsection{Algorithm: regularization with adaptive RKHS}
We estimate the vector $\bphi$ in \eqref{eq:discrete} by least squares, which minimizes the loss functional
\begin{align}\label{eq:lossFn_disc}
\widehat \bphi  = \argmin{\bphi\in \R^n}\costF(\bphi), \quad \quad  \costF(\bphi) &= \|\by- \bL \bphi\|^2 =  \sum_{i=1}^m |y_i- (\bL\bphi)_i |^2. 
\end{align}
With $\bA^{\dagger}$ denoting the pseudo-inverse of $\bA$, we compute a minimizer as
\begin{equation}\label{eq:lse_disc}
\widehat \bphi = \bA^{\dagger} \bb, \quad \text{ where }  \bA =  \bL^\top \bL, \, \bb =  \bL^\top \by .
\end{equation}

However, this inverse problem is ill-posed. The ill-posedness is rooted in the integral equation \eqref{eq:FIE}, as its solution involves the inversion of a compact operator with eigenvalues converging to zero \cite{nashed1974generalized}. Computationally, the ill-posedness is seen in the ill-conditionedness of the matrix $\bA$.

\begin{figure}[H]	\vspace{-2mm} 
    \centering 	
    {\includegraphics[width =0.45\textwidth]{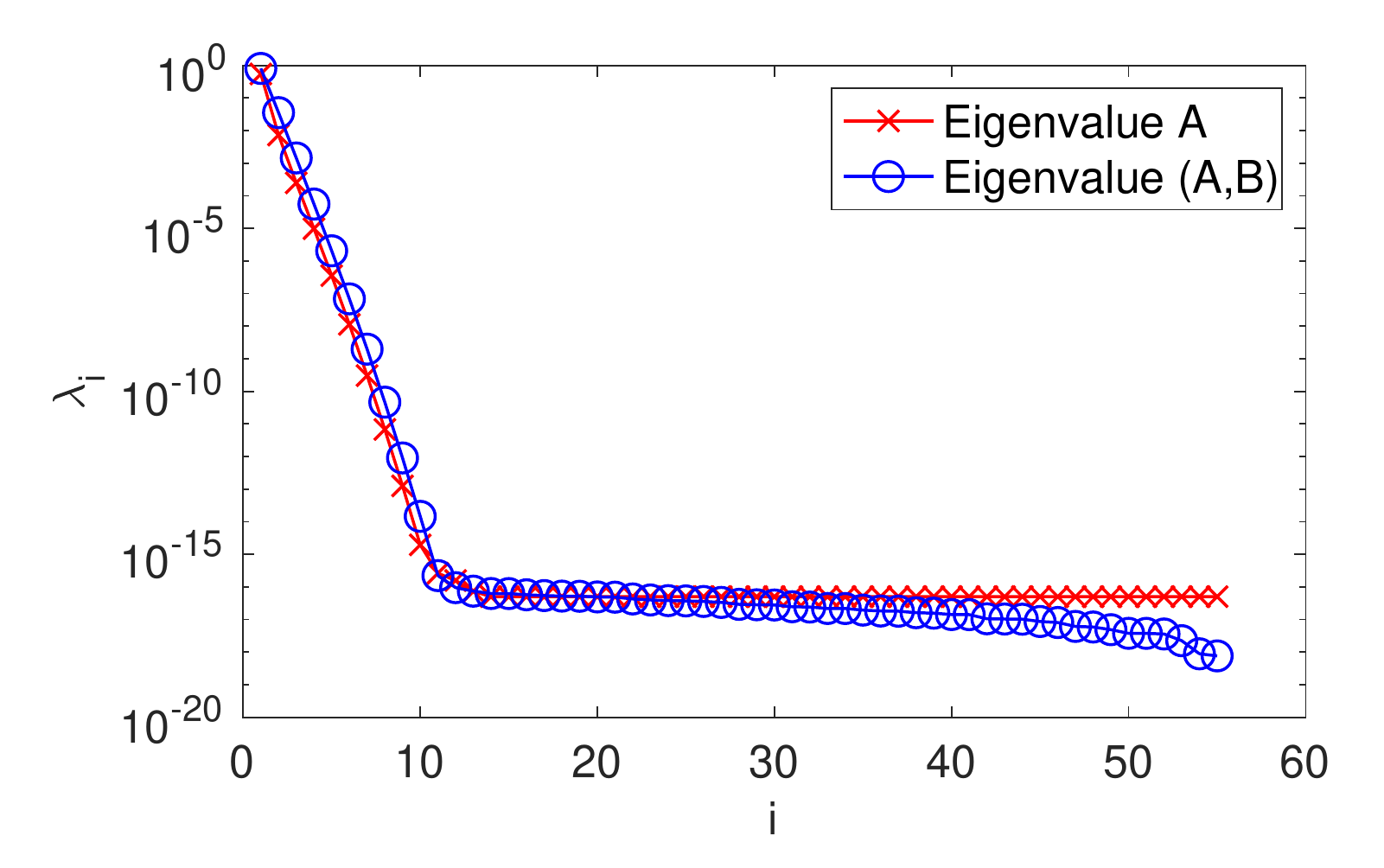}
    \includegraphics[width =0.45\textwidth]{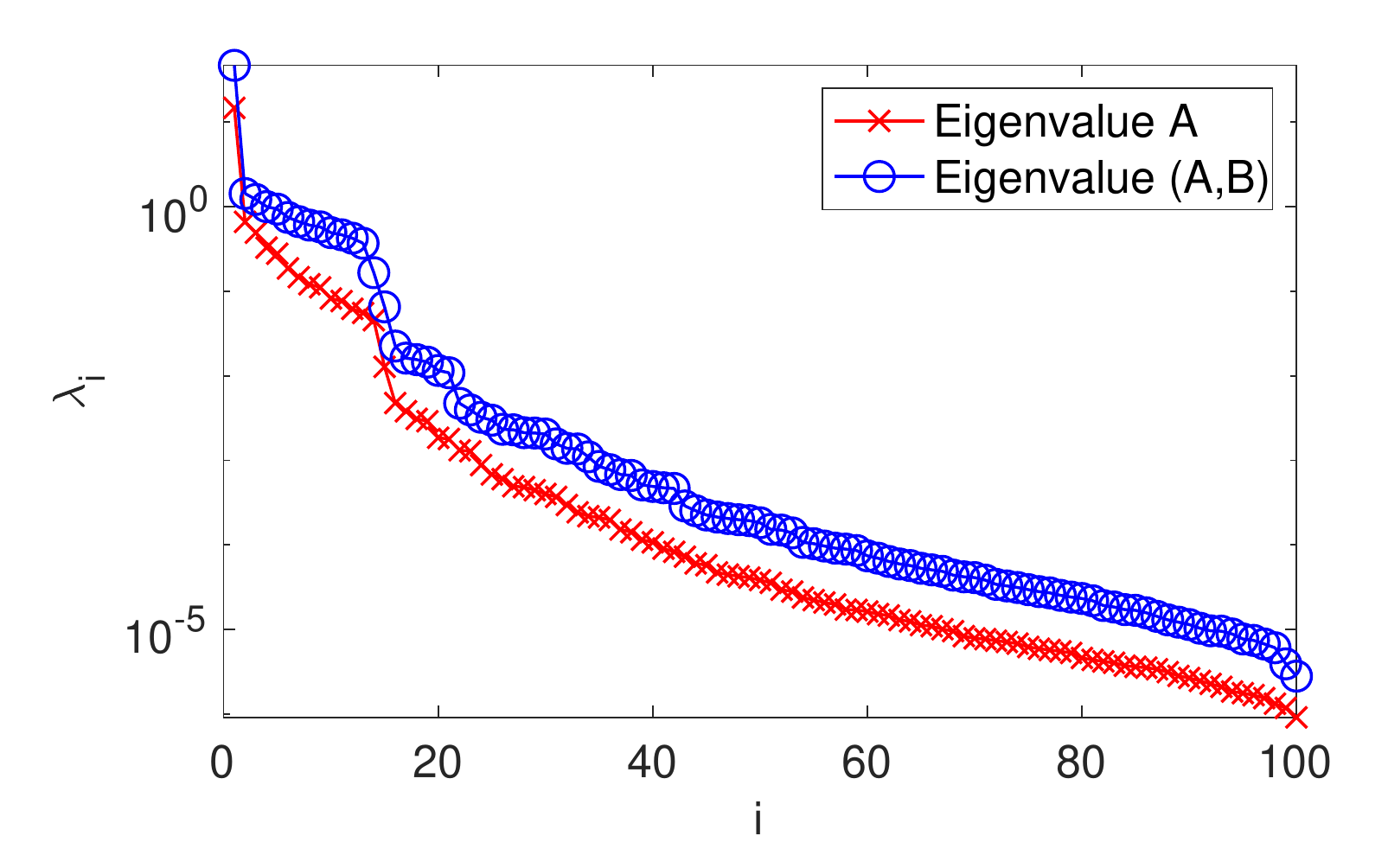} \\
    (a) Exponential decay \hspace{30mm} (b) Polynomial decay} 
\caption{Eigenvalues of $\bA$ and eigenvalues of $\LGbar$ corresponding to the two kernels in \eqref{eq:kernels}. {Here the eigenvalues of $\bA$ are computed by singular value decomposition and the eigenvalues of $\LGbar$ are computed as the generalized eigenvalues of $(\bA,\bB)$. }
} 
\label{fig:eigVals}		\vspace{-5mm}
\end{figure}

 Figure \ref{fig:eigVals} shows the eigenvalues of $\bA$ and the operator $\LGbar:L^2_\rho(\calS) \to L^2_\rho(\calS) $. Here the eigenvalues of $\bA$ are computed by singular value decomposition (SVD) to achieve better numerical stability than the eigenvalue decomposition. The eigenvalues of $\LGbar$ are computed by solving the generalized eigenvalue problem (see \cite[Theorem 4.1]{LLA22} or \cite[Proposition 5.6]{chada2022data})
\begin{equation}\label{eq:AbB1}
\bA V=  \bB\Lambda V, \quad s.t., V^\top \bB V = I_n, \quad \Lambda= \mathrm{Diag}(\lambda_1,\ldots,\lambda_n),
\end{equation}
which is computed by (symmetric) generalized eigenvalue decomposition \cite{nakatsukasa2013stable}. We observe that the generalized SVD method,
while producing stable generalized eigenvalues for these symmetric matrices, does not yield generalized eigenvectors satisfying the above equations. 
Hence, to take advantage of $\bA$ and $\bB$ being symmetric positive semi-definite, we use the generalized eigenvalue decomposition algorithm that yields eigenvectors satisfying these equations.    

\vspace{0.2in}

 Figure \ref{fig:eigVals} shows that the eigenvalues of $\bA$ and $\LGbar$ decay exponentially in $i$, providing an example of the assumption on $\lambda_i$ in Theorem \ref{thm:conv_sigma}. In particular, the fast-decaying eigenvalues show that the conditional numbers of $\bA$ and $\LGbar$ are extremely large, indicating severe ill-posedness. As a result, regularization is necessary to produce a stable solution from the noisy observations. 

We consider Tikhonov regularization that adds a quadratic penalty term ${\yolambda} \|\bphi\|_\bC^2 = {\yolambda} \bphi^\top \bC \bphi$ to the loss function, and solves the minimizer by least squares 
\begin{equation}\label{eq:lse_disc_reg}
\widehat \bphi_{\yolambda} = (\bA+{\yolambda} \bC)^{-1} \bb = \argmin{\bphi\in\R^n} \costF(\bphi) + {\yolambda} \|\bphi\|_\bC^2.
\end{equation} 
The hyper-parameter ${\yolambda}$ controls the strength of the regularization, and its selection has been thoroughly studied in \cite{wahba1977practical,hansen1998rank,hansen_LcurveIts_a}. We will select it by the L-curve method in \cite{hansen_LcurveIts_a}. 
The L-curve is a log-log plot of the curve 
$ l({\yolambda})=(y({\yolambda}),x({\yolambda}))$ with $y({\yolambda}) ^2=\| \bphi_{\yolambda}\|_\bC^2 $ and $x({\yolambda})^2 = \calE(\bphi_{\yolambda}) $ with $\bphi_{\yolambda}$ in \eqref{eq:lse_disc_reg}. The L-curve method maximizes the curvature to reach a balance between the minimization of the likelihood and the control of the regularization: 
\beqs 
	{\yolambda}_{*} 
	= {\rm argmax}_{\lambda_{\text{min}} \leq {\yolambda} \leq \lambda_{\text{max}}}\kappa(l ({\yolambda})),   \quad \kappa(l ({\yolambda}))= \frac{x'y'' - x' y''}{(x'\,^2 + y'\,^2)^{3/2}},
\eeqs 
where $\lambda_{min}$ and $\lambda_{max}$ are the smallest and the largest eigenvalues of $\bA$.

The main task is to select the norm $\|\cdot\|_\bC$. We consider the norm of the RKHS $H_G$ (which we explain in detail below), and we compare it with two commonly used norms: the Euclidean norm with $\bC=\mathbf{I}$ and the $L^2_\rho(\calS)$ norm with $\bC=\bB$ in \eqref{eq:basisMat}. We refer to them as the RKHS, $l^2$, and $L^2_\rho(\calS)$ regularization, respectively. They are summarized in Table \ref{tab:regularizers}. 
\begin{table}[htp] 
	\begin{center} 
		\caption{ Three regularizers using the norms of $l^2$, $L^2_\rho(\calS)$ and RKHS.   
		 } \label{tab:regularizers}
		\begin{tabular}{ l   l  l }		\toprule 
                 Regularizer name            & $\bC$  & Regularized estimator  \\  \hline
	         l2   &   $\mathbf{I}$  & $ \bphi_{\yolambda}^{l^2}= (\bA + {\yolambda} \mathbf{I})^{-1}\bb$ \\ 
	         L2    &   $\bB$  &  $ \bphi_{\yolambda}^{L^2}= (\bA + {\yolambda} \bB)^{-1}\bb$  \\
	       RKHS     &  $\bC_{rkhs}$  &  $ \bphi_{\yolambda}^{H_G}= (\bA + {\yolambda} \bC_{rkhs})^{-1}\bb$  \\
	 			\bottomrule	  
		\end{tabular}  \vspace{-4mm}
	\end{center}
\end{table}

\paragraph{The RKHS regularization.} Let $\{(\lambda_i,\bpsi_i)\}_{i=1}^n$ be the eigenvalue and eigen-function of $\LGbar$ over $L^2_\rho(\calS)$, and recall that $\{\bpsi_i\}$ form an orthonormal basis of $L^2_\rho(\calS)$. Here they are solved by the generalized eigenvalue problem in \eqref{eq:AbB1}, thus, $\bpsi_i = \sum_j V_{ji}\be_j$ and $\boldsymbol{\Psi} = (\bpsi_1,\ldots,\bpsi_n) = \mathbf{I} V$. Then, $\bphi= \mathbf{I} \bphi = \boldsymbol{\Psi} V^{-1} \bphi$. 
Thus, the $H_G$ norm of $\bphi $ is 
\begin{align*}
 \|\bphi\|_{H_G}^2 = &
\innerp{\bphi,\LGbar^{-1}\bphi}_{L^2_\rho(\calS)}  =  \innerp{\boldsymbol{\Psi} V^{-1} \bphi,\LGbar^{-1}\boldsymbol{\Psi} V^{-1} \bphi}_{L^2_\rho(\calS)}   \\
=  & ( V^{-1}\bphi)^\top \innerp{\boldsymbol{\Psi} , \LGbar^{-1}\boldsymbol{\Psi}}_{L^2_\rho(\calS)} V^{-1} \bphi = \bphi^\top  ( V^{-1})^\top \Lambda^{\dagger}V^{-1} \bphi \\
= & \|\bphi\|_{\bC_{rkhs}}^2 \, \text{ with } \bC_{rkhs} :=( V^{-1})^\top \Lambda^{\dagger}V^{-1}.
 \end{align*}
In particular, if $\rho$ is a uniform measure and the basis matrix is $\bB=I_n$, we have $\bC_{rkhs}= \bA^{\dagger}$.

In computation, we improve numerical stability by avoiding the pseudo-inverse of a singular matrix. The procedure is as follows. Note that for $\bC_* := V \Lambda^{1/2} $, we have $\bC_*^\top  \bC_{rkhs} \bC_* = \begin{pmatrix}I_r & 0 \\ 0 & 0  \end{pmatrix} := \mathbf{I}_r$, where $I_r$ is the identity matrix with rank $r$, the number of positive eigenvalues in $\Lambda$.  
Then, the linear equation $(\bA + {\yolambda} \bC_{rkhs}) c_{\yolambda} =\bb$, the equation \eqref{eq:lse_disc_reg} with regularization matrix $\bC_{rkhs}$, is equivalent to 
\begin{equation}\label{eq:lsqminnorm}
(\bC_*\bA \bC_*+ {\yolambda} \mathbf{I}_r) \widetilde \bphi_{\yolambda} = \bC_*\bb
\end{equation}
 with $\widetilde \bphi_{\yolambda}= \bC_*^{-1}\bphi_{\yolambda}$. We compute $\widetilde \bphi_{\yolambda}$ in the above equation by least squares with the minimal norm; then we obtain $\bphi_{\yolambda}= \bC_* \widetilde \bphi_{\yolambda}$. These treatments avoid the inversions of ill-conditioned or singular matrices and lead to more robust estimators.

The implementation of the RKHS regularization based on matrix decomposition is summarized in Algorithm \ref{alg:dartr}. 
\begin{algorithm}[H]
{\small
\begin{algorithmic}[1]
\Require{The regression triplet $(\bA, \bb, \bB)$ consisting of normal matrix $\bA$, vector $\bb$ in \eqref{eq:lse_disc} and basis matrix $\bB$ as in \eqref{eq:basisMat}.} 
\Ensure{An RKHS regularized estimator $\widehat \bphi_{{\yolambda}_0}$ and its loss value $\calE(\widehat \bphi_{{\yolambda}_0})$. }
	 \State Solve the generalized eigenvalue problem $\bA V = \bB V\Lambda$, $V^\top \bB V = I$.  
	 \State Compute $\widetilde \bA= \bC_{*}^{\top}\bA \bC_{*}$ and $\widetilde \bb = \bC_{*} \bb $ with $ \bC_{*} = V \Lambda^{1/2}$. 	 
	\State Use the L-curve method to find an optimal estimator $\widehat \bphi_{{\yolambda}_0}$: 
	\begin{itemize}
	\item Set the range for ${\yolambda}$ to be the range of the eigenvalues in $\Lambda$. 
	\item For each ${\yolambda}$,  solve $\widetilde \bphi_{\yolambda}$ from $(\widetilde \bA +{\yolambda} I) \widetilde \bphi_{\yolambda} = \widetilde\bb$ by least squares with minimal norm and set $ \widehat \bphi_{\yolambda} = \bC_* \widetilde \bphi_{\yolambda} $. 
	\item  Select ${\yolambda}_0$ maximizing the curvature of the ${\yolambda}$-curve $(\log \calE( \widehat \bphi_{\yolambda}), \log(\widehat \bphi_{\yolambda}^\top\bC_{rkhs} \widehat \bphi_{\yolambda} ))$.
	\item Return the solution $\widetilde \bphi_{{\yolambda}_0}$ solving $(\widetilde \bA +{\yolambda}_0 I) \widetilde \bphi_{{\yolambda}_0} = \widetilde\bb$. 
	\end{itemize}
\end{algorithmic}
\caption{Solving $\bA \bphi\approx \bb$ with basis matrix $\bB$ by the adaptive RKHS Regularization. 
 }\label{alg:dartr}
}
\end{algorithm}\vspace{-2mm}

\paragraph{Function space of identifiability (FSOI).} Recall that the FSOI $H$ 
is the linear subspace of $L^2_\rho(\calS)$ spanned by the eigenvectors of $\LGbar$ with nonzero eigenvalues. 
In computation, the FSOI is the eigenspace of $\LGbar$ spanned by the eigenvectors with positive eigenvalues (greater than the numerical precision in practice). Figure \ref{fig:eigVals} suggests that the FSOI is low-dimensional and contains only a few eigenvectors.

Also, we access the performance of an estimator by its accuracy in estimating the projection of the true solution in the FSOI, 
since the data provides no information for estimating the components outside the FSOI. That is, we compute the error of an estimator $\bphi$ as 
\begin{equation}\label{eq:error}
Err = \|  P_H \widehat{\bphi} - P_H \bphi_* \|_{L^2_\rho(\calS)}^2, 
\end{equation} 
where $P_H$ is the projection to $H$.

\begin{figure}[t]	\vspace{-2mm} 
    \centering 	
    \includegraphics[width =0.99\textwidth]{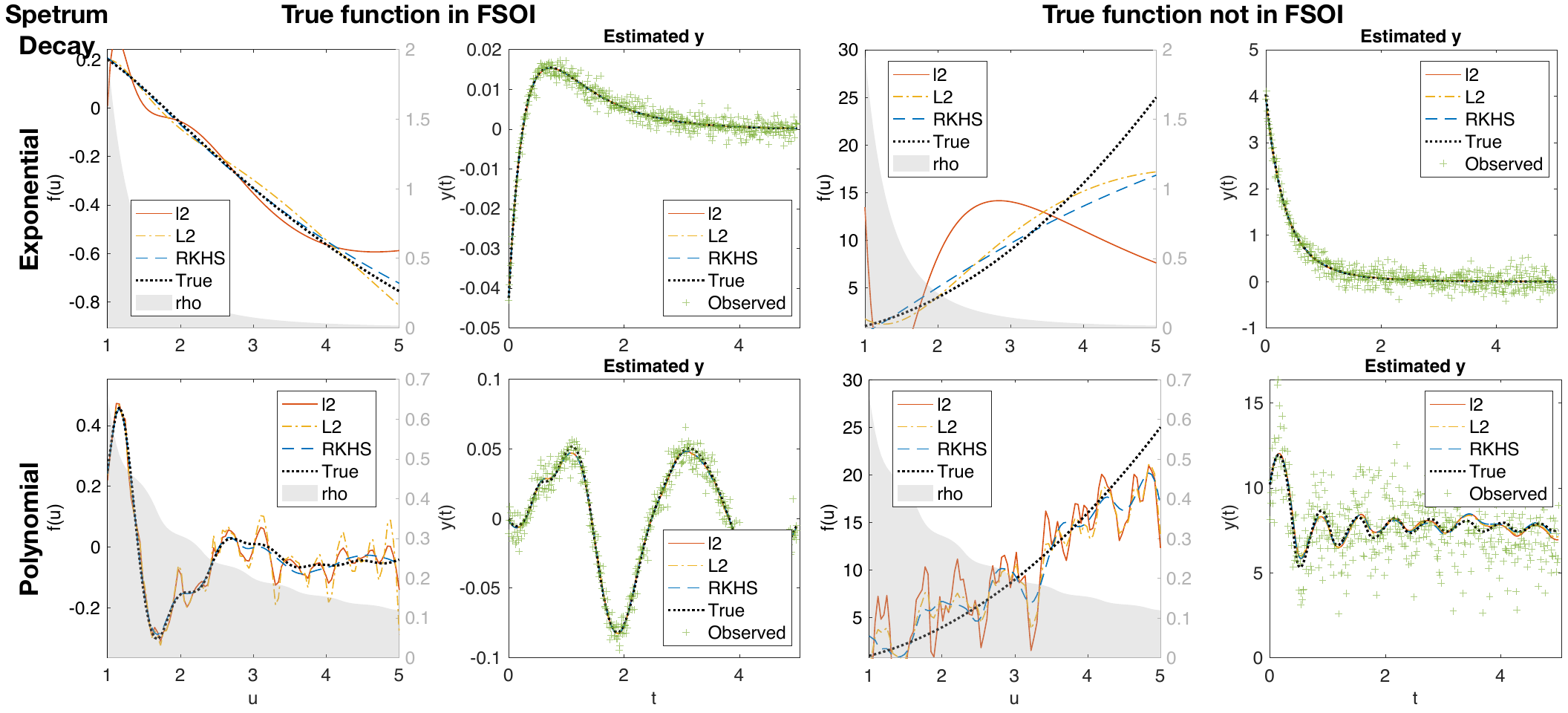}
	\vspace{-2mm} 
\caption{Typical regularized estimators when $nsr=2$ and their recovery of the signal. 
All estimators recover the true signal accurately. 
 The RKHS regularizer significantly outperforms the other two regularizers when the true function is in the FSOI (see top-left and bottom-left), but it may slightly underperform the $L^2_\rho(\calS)$ regularizer when the true function is not in the FSOI (see top- and bottom-right). 
} 
\label{fig:estimator}		\vspace{-5mm}
\end{figure}

\subsection{Convergence as the noise level decays}\label{sec:conv_noise}
To test the convergence of the regularized estimators as the noise level decays, we set the noise-to-signal ratio to be $nsr \in \{0.125,0.25,0.5,1,2\}$ and $\Delta t = 0.01$. Recall that $\sigma= \|\bL \bphi  \|  \times nsr$.

We consider both the scenario of $\phi_{true}\in \mbox{FSOI}$ and the scenario of $\phi_{true}\notin \mbox{FSOI}$ {or it has significant components in the eigenspace of small eigenvalues. We call these scenarios ``in FSOI'' and ``not in FSOI'' for simplicity.} For the first scenario, we set $\phi_{true}=\psi_2$, the second eigenvector of $\LGbar$. For the second scenario, we consider $\phi_{true}=x^2$. We note that the latter is challenging to recover because {of $\phi_{true}$'s components in eigenspace with small or zero eigenvalues.}

Figure \ref{fig:estimator} shows the typical regularized estimators in the two scenarios and their de-noised output $\widehat{y} = \bL \widehat \bphi$. The left column shows that when the true solution is in the FSOI, the RKHS regularizer significantly outperforms the other two.  However, all the regularized estimators can recover the true signal (or de-noise the data) accurately, even though these estimators may deviate largely from the true solution. The right column shows that when the true solution has components outside the FSOI, the RKHS regularizer can be less accurate than the $L^2_\rho(\calS)$ regularizer. At the same time, both significantly outperform the $l^2$ regularizer. Yet again, all regularized estimators can de-noise the data accurately even when they have large errors. Thus, this inverse problem is ill-defined, and one must restrict the inverse to be in the FSOI.

\begin{figure}[htb]	\vspace{-2mm} 
    \centering 	
    {\includegraphics[width =0.98\textwidth]{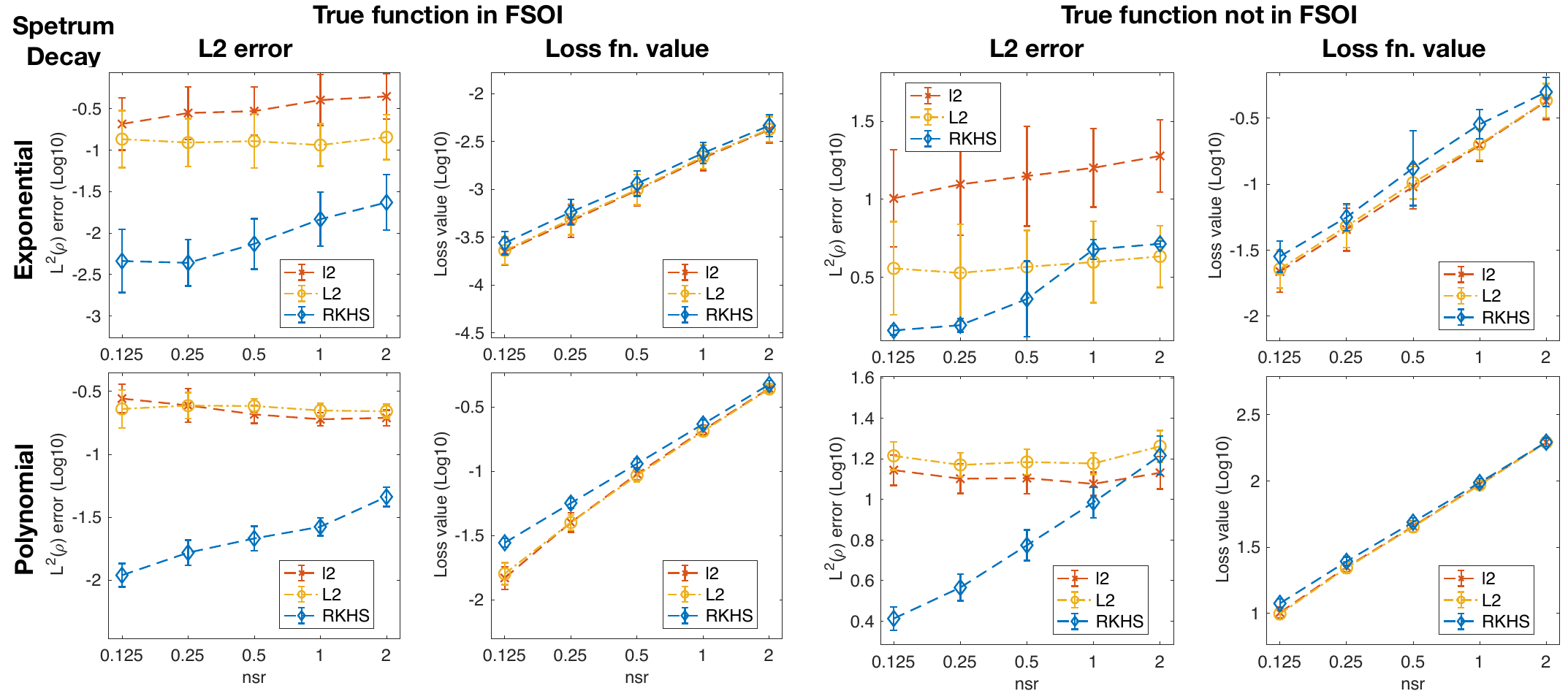}}
	\vspace{-2mm}  
\caption{Mean and standard deviations of the projected error $Err$ in \eqref{eq:error} and the loss functional values in 100 simulations. The RKHS regularizer significantly outperforms the other two regularizers for $\phi_{true} \in \mbox{FSOI}$ (see top- and bottom-left). In all cases, the RKHS regularized estimator has errors consistently decaying as the noise level decreases, and it yields minimal loss functional values that are slightly larger than those of the $l^2$ and $L^2$-regularized estimators.   
 } \label{fig:conv_nsr}		\vspace{-3mm}
\end{figure}

Figure \ref{fig:conv_nsr} demonstrates the convergence of the regularized estimators and their loss values as the noise-to-signal-ratio decreases, with $nsr \in \{0.125,0.25,0.5,1,2\}$. It presents the mean and standard deviations of the regularized estimators' $L^2_\rho(\calS)$ errors and the loss values {in 100 simulations, each with a different noise realization}. The errors are computed for the estimators' projection inside the FSOI as in \eqref{eq:error}.  When the true solution is inside the FSOI, the RKHS regularized estimator has errors significantly smaller than those of the other two regularizers, though it has the largest loss values. When the true solution is outside the FSOI, the RKHS regularized estimator has slightly smaller errors than the $L^2_\rho(\calS)$ regularizer, and both outperform the $l^2$ regularizer. Importantly, when the true solution is inside FOSI, the RKHS-regularized estimator shows a clear linear error decay in noise for $nsr\in \{0.25,0.5,1,2\}$. In contrast, the $L^2_\rho$-regularized estimator has an error remaining at a constant level.

\begin{figure}[htb]	\vspace{-2mm} 
    \centering 	
    {\includegraphics[width =0.98\textwidth]{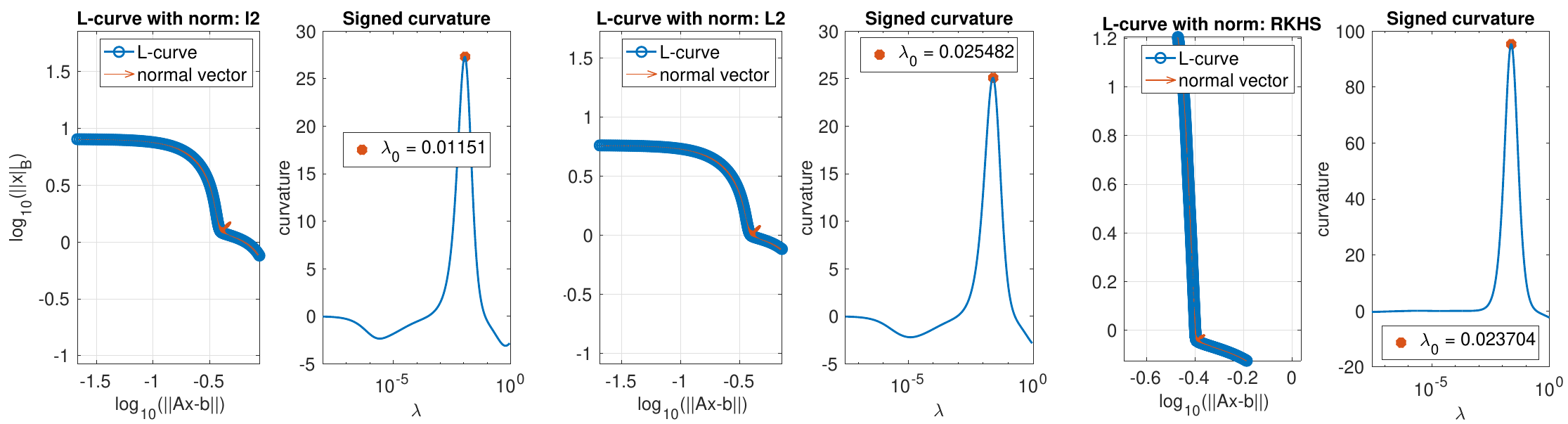}}
	\vspace{-2mm}  
\caption{L-curves in a typical simulation for the $l^2$-, $L^2$- and RKHS-regularizers. The RKHS-regularizer has a better L-shaped curve than the other two regularizers, indicating that the RKHS-norm is a better metric for regularization than the other two norms. 
 } \label{fig:Lcurve}		\vspace{-1mm}
\end{figure}

Notably in Figure \ref{fig:conv_nsr}, the RKHS-regularizer achieves a rate of convergence about $\sigma^{0.5}$ for the $Err$ defined in \eqref{eq:error}, and this rate agrees with the theoretical rate for the squared error in Theorem \ref{thm:conv_sigma}. Such a close agreement between theory and practical computation is amazing since the theoretical rate uses the oracle optimal hyper-parameter based on the true solution, but the practical computation empirically selects the hyper-parameter by the L-curve method. In contrast, the $L^2$-regularizer does not yield a convergent estimator. Since the two regularizers only differ at the choice of regularization norm $\| \phi \|_\square$, we conclude that the RKHS-norm provides a better metric for the L-curve than the $L^2_\rho$ norm. Figure \ref{fig:Lcurve} shows that the L-curve in the RKHS-norm is in a better shape than the other two norms.     

\subsection{Convergence as the observation mesh refines}
Next, we test the convergence of the regularized estimator as the observation mesh refines.  We set $\Delta t\in 0.005\times \{1,2,4,8,16\}$ and $nsr =1$. 

We continue to consider the two scenarios as in Section \ref{sec:conv_noise} with the same choice of the true functions. However, we compute the FSOI using a finer observation mesh $\Delta t = 0.0005$.

\begin{figure}[H]	\vspace{-2mm} 
    \centering 	
    \includegraphics[width =0.98\textwidth]{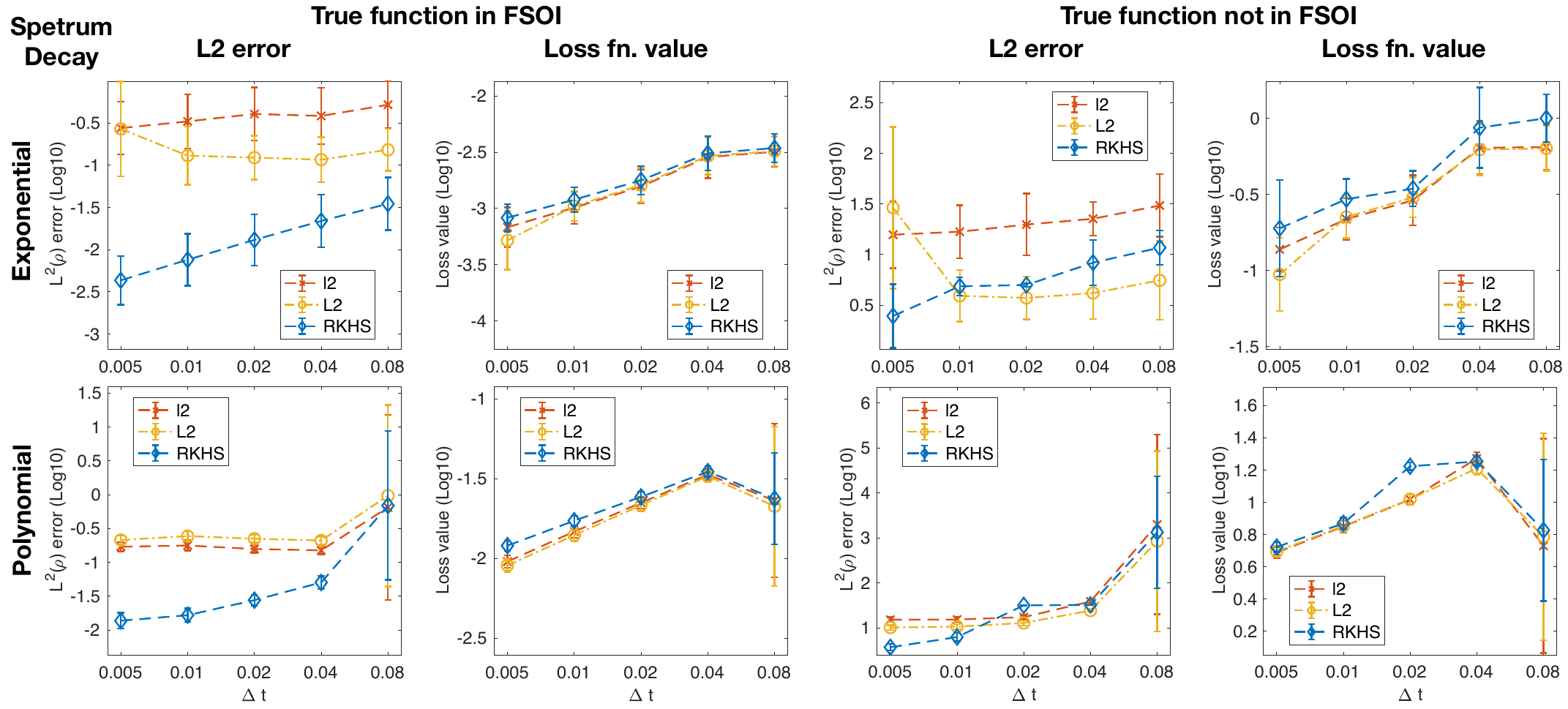}
	\vspace{2mm}  
\caption{Mean and standard deviations of the projected error $Err$ in \eqref{eq:error} and the loss functional values in 100 independent simulations. 
 The RKHS regularizer significantly outperforms the other two regularizers when $\phi_{true}\in \mbox{FSOI}$ (see top-and bottom-left). In all cases, only the RKHS regularized estimator has errors decaying consistently with the decrease of $\Delta t$. 
 } \label{fig:conv_dt}		\vspace{-5mm}
\end{figure}

Figure \ref{fig:conv_dt} demonstrates the convergence of the regularized estimators and their loss values as the observation mesh increases. It presents the mean and standard deviations of the projected error $Err$ defined in \eqref{eq:error} and the loss values {in 100 simulations, each with an independent noise realization}. The RKHS regularized estimator has errors consistently decaying as the mesh refines, outperforming the other two regularizers when the true solution is inside the FSOI (top- and bottom-left), although its loss functional values are slightly larger. Also, the RKHS-regularizer slightly underperforms the $L^2$-regularizer when the true solution is outside the FSOI, and both outperforms the $l^2$ regularizer. Importantly, in all the cases, the RKHS regularized estimator has errors consistently decaying as the time mesh refines.

\section{Discussions and limitations}     \label{sec:limitations}

\paragraph{Comparative analysis.} A major contribution of this study is the comparative analysis between the $L^2$-regularized estimator and the $H_G$-regularized estimator via the sharp convergence rate in the small noise limit. The sharp convergence rate of the expected $L^2_\rho$ error depends on two factors: spectrum decay of the operator of inversion, and the speed of decay of the true function's cofficients. The key innovation is extracting these factors from multiple mixed factors in regularization. In particular, the expected error under the Gaussian noise assumption helps remove the randomness in measurement error, and the oracle optimal hyper-parameter minimizing the expected $L^2_\rho$ error removes the uncertainty associated with the empirical selection of the hyper-parameter.    

 Our comparative analysis can be extended in two directions.  
        \begin{itemize}
        \item  \emph{Comparing other Hilbert space norms for regularization.} In this study, we compared the RKHS-regularizer with the widely-used $L^2$-regularizer. But our small noise analysis can also be applied to comparing other regularization norms, such as other RKHS norms with pre-selected reproducing kernels or the $H^1$ norm.
        However, it may not be suitable to analyze non-Hilbert space norms, since the analysis presented here relies on the spectral properties of the self-adjoint operator of inversion. 
        \item \emph{Convergence with respect to mesh refinement}. It remains open to studying the convergence of the regularized estimator when the observation mesh refines. The challenge comes from a complicated balance between three elements: the discretization error, the observation noise-induced error, and the regularization-induced bias.
\end{itemize}
         
\paragraph{Practical computations.} Another contribution of this study is the new adaptive RKHS regularization for solving the Fredholm integral equation of the first kind. This study implements the RKHS regularization by a matrix decomposition-based algorithm. Numerical results show that the RKHS-regularizer outperforms the widely used $l^2$- and $L^2$-regularizers. However, the matrix decomposition-based method will be computationally costly when the size of the matrix is large. For large-scale problems, iterative methods are necessary, and the challenge comes from the computation of the RKHS norm and a proper early stopping strategy. These investigations are beyond the scope of this study, and we leave it for future work. 


\section{Conclusion and future work}\label{sec:conlusion}
We have introduced an adaptive RKHS regularization for ill-posed linear inverse problems. The adaptive RKHS arises in the variational formulation of the inverse problem, and it is determined by the linear operator and the observation mesh. Its closure is the space in which we can identify the true solution. 

Importantly, we have introduced a quantitative comparative analysis to compare different regularization norms through the sharp convergence rate of the mean-square error of the regularized estimator in the small noise limit. We have proved that the RKHS-regularized and the widely-used $L^2$-regularized estimators have the same sharp convergence rate, but the RKHS-regularized estimator has a smaller multiplicative factor. 

Additionally, we present an implementation of the RKHS-regularization by an algorithm based on matrix decomposition. The algorithm uses the L-curve method to select the optimal hyperparameter. Numerical tests confirm the robust convergence of the RKHS-regularized estimator when either the noise level decays or the observation mesh refines, no matter whether the true function is in the FSOI or not. In the tests, the RKHS regularizer significantly outperforms the widely used $l^2$- and $L^2$-regularizers. This implementation relies on the matrix decomposition, which {is costly and} sensitive to numerical error when {the matrix size is large}. A future direction is to develop iterative methods that minimize the penalized loss functional via iterations with an early stopping strategy \cite{hansen1998rank,arioli2013GKb}, which is suitable for solving {large sized} problems.

 \paragraph{Acknowledgments} The work of FL is partially funded by the Johns Hopkins University Catalyst Award and AFOSR FA9550-21-1-0317 and FA9550-20-1-0288. FL would like to thank Quanjun Lang for helpful discussions on the L-curve methods and the convergence tests. The work of MYO is partially supported by the European Union's Horizon 2020 Research and Innovation Programme under the Marie Sk\l odowska-Curie grant agreement 101008231 (EXPOWER).
 
 \bibliographystyle{plain}

\begin{thebibliography}{10}

\bibitem{arioli2013GKb}
Mario Arioli.
\newblock Generalized golub--kahan bidiagonalization and stopping criteria.
\newblock {\em SIAM Journal on Matrix Analysis and Applications},
  34(2):571--592, 2013.

\bibitem{bauer2007regularization}
Frank Bauer, Sergei Pereverzev, and Lorenzo Rosasco.
\newblock On regularization algorithms in learning theory.
\newblock {\em Journal of complexity}, 23(1):52--72, 2007.

\bibitem{Bi2022span-of-regular}
Chuan Bi, M.~Yvonne Ou, Mustapha Bouhrara, and Richard~G. Spencer.
\newblock Span of regularization for solution of inverse problems with
  application to magnetic resonance relaxometry of the brain.
\newblock {\em Scientific Reports}, 12(1):20194, 2022.

\bibitem{chada2022data}
Neil~K Chada, Quanjun Lang, Fei Lu, and Xiong Wang.
\newblock A data-adaptive prior for {Bayesian} learning of kernels in
  operators.
\newblock {\em arXiv preprint arXiv:2212.14163}, 2022.

\bibitem{chen2022stochastic}
Zhiming Chen, Wenlong Zhang, and Jun Zou.
\newblock Stochastic convergence of regularized solutions and their finite
  element approximations to inverse source problems.
\newblock {\em SIAM Journal on Numerical Analysis}, 60(2):751--780, 2022.

\bibitem{cucker2007learning-theory}
Felipe Cucker and Ding~Xuan Zhou.
\newblock {\em Learning Theory: An Application Theory Viewpoint}.
\newblock Cambridge Monographs on Applied and Computational Methematics.
  Cambridge University Press, 1 edition, 2007.

\bibitem{CZ07book}
Felipe Cucker and Ding~Xuan Zhou.
\newblock {\em Learning theory: an approximation theory viewpoint}, volume~24.
\newblock Cambridge University Press, {Cambridge}, 2007.

\bibitem{engl1996regularization}
Heinz~Werner Engl, Martin Hanke, and Andreas Neubauer.
\newblock {\em Regularization of inverse problems}, volume 375.
\newblock Springer Science \& Business Media, 1996.

\bibitem{gazzola2019ir}
Silvia Gazzola, Per~Christian Hansen, and James~G Nagy.
\newblock Ir tools: a matlab package of iterative regularization methods and
  large-scale test problems.
\newblock {\em Numerical Algorithms}, 81(3):773--811, 2019.

\bibitem{hadamard1923lectures}
Jacques~Salomon Hadamard.
\newblock {\em Lectures on Cauchy's problem in linear partial differential
  equations}, volume~18.
\newblock Yale University Press, 1923.

\bibitem{hansen1998rank}
Per~Christian Hansen.
\newblock {\em Rank-deficient and discrete ill-posed problems: numerical
  aspects of linear inversion}.
\newblock SIAM, 1998.

\bibitem{hansen_LcurveIts_a}
Per~Christian Hansen.
\newblock The {L}-curve and its use in the numerical treatment of inverse
  problems.
\newblock In {\em in Computational Inverse Problems in Electrocardiology, ed.
  P. Johnston, Advances in Computational Bioengineering}, pages 119--142. WIT
  Press, 2000.

\bibitem{LangLu21}
Quanjun Lang and Fei Lu.
\newblock Identifiability of interaction kernels in mean-field equations of
  interacting particles.
\newblock {\em Foundations of Data Science}, 2023.

\bibitem{LangLu23sna}
Quanjun Lang and Fei Lu.
\newblock Small noise analysis for tikhonov and rkhs regularizations.
\newblock {\em arXiv preprint arXiv:2305.11055}, 2023.

\bibitem{Li2005modified}
Gongsheng Li and Zuhair Nashed.
\newblock A modified {Tikhonov} regularization for linear operator equations.
\newblock {\em Numerical functional analysis and optimization},
  26(4-5):543--563, 2005.

\bibitem{LAY22}
Fei Lu, Qingci An, and Yue Yu.
\newblock Nonparametric learning of kernels in nonlocal operators.
\newblock {\em J. Peridyn Nonlocal Model}, 2023.

\bibitem{LLA22}
Fei Lu, Quanjun Lang, and Qingci An.
\newblock {Data adaptive RKHS Tikhonov regularization for learning kernels in
  operators}.
\newblock {\em Proceedings of Mathematical and Scientific Machine Learning,
  PMLR 190:158-172}, 2022.

\bibitem{nakatsukasa2013stable}
Yuji Nakatsukasa and Nicholas~J Higham.
\newblock Stable and efficient spectral divide and conquer algorithms for the
  symmetric eigenvalue decomposition and the svd.
\newblock {\em SIAM Journal on Scientific Computing}, 35(3):A1325--A1349, 2013.

\bibitem{nashed1974generalized}
M~Zuhair Nashed and Grace Wahba.
\newblock Generalized inverses in reproducing kernel spaces: an approach to
  regularization of linear operator equations.
\newblock {\em SIAM Journal on Mathematical Analysis}, 5(6):974--987, 1974.

\bibitem{rudin1992nonlinear}
Leonid~I Rudin, Stanley Osher, and Emad Fatemi.
\newblock Nonlinear total variation based noise removal algorithms.
\newblock {\em Physica D: nonlinear phenomena}, 60(1-4):259--268, 1992.

\bibitem{tibshirani1996_RegressionShrinkage}
Robert Tibshirani.
\newblock Regression shrinkage and selection via the {{Lasso}}.
\newblock {\em Journal of the Royal Statistical Society: Series B
  (Methodological)}, 58(1):267--288, 1996.

\bibitem{tihonov1963solution}
Andrei~Nikolajevits Tikhonov.
\newblock Solution of incorrectly formulated problems and the regularization
  method.
\newblock {\em Soviet Math.}, 4:1035--1038, 1963.

\bibitem{wahba1973convergence}
Grace Wahba.
\newblock Convergence rates of certain approximate solutions to fredholm
  integral equations of the first kind.
\newblock {\em Journal of Approximation Theory}, 7(2):167--185, 1973.

\bibitem{wahba1977practical}
Grace Wahba.
\newblock Practical approximate solutions to linear operator equations when the
  data are noisy.
\newblock {\em SIAM journal on numerical analysis}, 14(4):651--667, 1977.

\bibitem{zhang2022stochastic}
Ye~Zhang and Chuchu Chen.
\newblock Stochastic asymptotical regularization for linear inverse problems.
\newblock {\em Inverse Problems}, 39(1):015007, 2022.

\end{thebibliography}

\end{document}